\documentclass[12pt]{amsart}
\usepackage{amsmath,amscd,amssymb,amsfonts}
\newcommand{\h}{\hbox}
\newcommand{\q}{\quad}
\newcommand{\nin}{\par\noindent}
\newcommand{\bs}{\par\bigskip}
\newcommand{\ms}{\par\medskip}
\newcommand{\sk}{\par\smallskip}
\newcommand{\bsn}{\par\bigskip\noindent}
\newcommand{\msn}{\par\medskip\noindent}
\newcommand{\skn}{\par\smallskip\noindent}
\newcommand{\bl}{\bigl}
\newcommand{\br}{\bigl}
\newcommand{\ssb}{\raise.15ex\h{${\scriptscriptstyle\bullet}$}}
\newcommand{\ssc}{\,\raise.15ex\h{${\scriptstyle\circ}$}\,}
\newcommand{\msum}{\h{$\sum$}}
\newcommand{\mopl}{\h{$\bigoplus$}}
\newcommand{\mcap}{\h{$\bigcap$}}
\newcommand{\mcup}{\h{$\bigcup$}}
\newcommand{\mprod}{\h{$\prod$}}

\newcommand{\ab}{{\bf a}}
\newcommand{\al}{\alpha}

\newcommand{\C}{{\mathbb C}}
\newcommand{\dd}{\partial}
\newcommand{\ddd}{{\rm d}}
\newcommand{\e}{{\mathbf e}}
\newcommand{\B}{{\mathcal B}}
\newcommand{\BB}{\widetilde{\mathcal B}}
\newcommand{\D}{{\mathcal D}}
\newcommand{\DD}{{\mathbb D}}
\newcommand{\ep}{\varepsilon}
\newcommand{\F}{{\mathcal F}}
\newcommand{\G}{{\mathcal G}}
\newcommand{\HH}{{\mathbf H}}
\newcommand{\Hc}{{\mathcal H}}
\newcommand{\Ht}{\widetilde{H}}
\newcommand{\I}{{\mathcal I}}
\newcommand{\J}{{\mathcal J}}
\newcommand{\JC}{{\rm JC}}

\newcommand{\la}{\lambda}
\newcommand{\La}{\Lambda}

\newcommand{\lct}{{\rm lct}}
\newcommand{\M}{{\mathcal M}}
\newcommand{\MM}{{\mathcal M}^{\ssb}}
\newcommand{\N}{{\mathbb N}}
\newcommand{\OO}{{\mathcal O}}
\newcommand{\PP}{{\mathbb P}}
\newcommand{\Q}{{\mathbb Q}}

\newcommand{\R}{{\mathbb R}}
\newcommand{\ft}{\widetilde{f}}

\newcommand{\Sc}{{\mathcal S}}
\newcommand{\Sci}{S^{\circ}}
\newcommand{\Si}{\Sigma}

\newcommand{\Uc}{{\mathcal U}}

\newcommand{\X}{\widetilde{X}}
\newcommand{\Yt}{\widetilde{Y}}
\newcommand{\y}{\widetilde{y}}
\newcommand{\Xc}{{\mathcal X}}
\newcommand{\Xcc}{{\mathcal X}^{\circ}}
\newcommand{\Ycc}{{\mathcal Y}^{\circ}}
\newcommand{\om}{\omega}
\newcommand{\omt}{\widetilde{\omega}}
\newcommand{\Z}{{\mathbb Z}}
\newcommand{\can}{{\rm can}}
\newcommand{\EM}{\widetilde{\rm EM}}
\newcommand{\BM}{{\rm BM}}
\newcommand{\CH}{{\rm CH}}
\newcommand{\DR}{{\rm DR}}
\newcommand{\Gr}{{\rm Gr}}
\newcommand{\MHM}{{\rm MHM}}
\newcommand{\Sing}{{\rm Sing}}
\newcommand{\Supp}{{\rm Supp}}
\newcommand{\spr}{{\rm sp}}
\newcommand{\vir}{{\rm vir}}

\newcommand{\into}{\hookrightarrow}
\newcommand{\simto}{\buildrel\sim\over\longrightarrow}
\newcommand{\ges}{\geqslant}
\newcommand{\les}{\leqslant}
\begin{document}
\title[Spectral Hirzebruch-Milnor classes]
{Spectral Hirzebruch-Milnor classes\\of singular hypersurfaces}
\author[L. Maxim ]{Laurentiu Maxim}
\address{L. Maxim : Department of Mathematics, University of Wisconsin-Madison, 480 Lincoln Drive, Madison WI 53706-1388 USA}
\email{maxim@math.wisc.edu}
\author[M. Saito ]{Morihiko Saito}
\address{M. Saito: RIMS Kyoto University, Kyoto 606-8502 Japan}
\email{msaito@kurims.kyoto-u.ac.jp}
\author[J. Sch\"urmann ]{J\"org Sch\"urmann}
\address{J. Sch\"urmann : Mathematische Institut, Universit\"at M\"unster, Einsteinstr. 62, 48149 M\"unster, Germany}
\email{jschuerm@uni-muenster.de}
\begin{abstract} We introduce spectral Hirzebruch-Milnor classes for singular hypersurfaces. These can be identified with Steenbrink spectra in the isolated singularity case, and may be viewed as their global analogues in general. Their definition uses vanishing cycles of mixed Hodge modules and the Todd class transformation. These are compatible with the pushforward by proper morphisms, and the classes can be calculated by using resolutions of singularities. Formulas for Hirzebruch-Milnor classes of projective hypersurfaces in terms of these classes are given in the case where the multiplicity of a generic hyperplane section is not 1. These formulas using hyperplane sections instead of hypersurface ones are easier to calculate in certain cases. Here we use the Thom-Sebastiani theorem for the underlying filtered $D$-modules of vanishing cycles, from which we can deduce the Thom-Sebastiani type theorem for spectral Hirzebruch-Milnor classes. For the Chern classes after specializing to $y=-1$, we can give a relatively simple formula for the localized Milnor classes, which implies a new formula for the Euler numbers of projective hypersurfaces, using iterated hyperplane sections. Applications to log canonical thresholds and Du Bois singularities are also explained; for instance, the latter can be detected by using Hirzebruch-Milnor classes in the projective hypersurface case.
\end{abstract}
\maketitle
\centerline{\bf Introduction}
\bs\nin
Let $Y$ be a smooth complex projective variety with $L$ a very ample line bundle on $Y$. Let $X$ be a hypersurface section of $Y$ defined by $s\in\Gamma(Y,L^{\otimes m})\setminus\{0\}$ for some $m\in\Z_{>0}$.
When $m=1$, a formula for the Hirzebruch-Milnor class (which expresses the difference between the Hirzebruch class and the virtual one) was given in \cite{MSS1} by using a sufficiently general section of $L$. Specializing to $y=-1$ and using \cite[Proposition 5.21]{Sch3}, this implies a formula for the Chern-Milnor class conjectured by S.~Yokura \cite{Yo2}, and proved by A.~Parusi\'nski and P.~Pragacz \cite{PP} (where $m=1$).
\sk
It is, however, desirable to generalize the above formula to the case $m>1$, since it is sometimes easier to calculate hyperplane sections than hypersurface ones, see Corollary~1 and Example~(2.8) below. In order to realize this, we need an inductive argument as below.
\sk
We denote the Hirzebruch class and the virtual one by $T_{y*}(X)$, $T^{\,\vir}_{y*}(X)\in\HH_{\ssb}(X)[y]$ for $X$ as above (see also \cite{MSS1}), where $\HH_k(X)=H_{2k}^{\BM}(X,\Q)$ or $\CH_k(X)_{\Q}$, see (1.1) below.
For $g\in\Gamma(Y,\OO_Y)$ with $Y$ a smooth complex algebraic variety in general, we denote by $\varphi_g\Q_{h,Y}$ the mixed Hodge module on ${\rm Sing}\,g^{-1}(0)$ up to a shift of complex such that its underlying $\Q$-complex is the vanishing cycle complex $\varphi_g\Q_Y$ in the sense of \cite{De2}, see \cite{mhp}, \cite{mhm}.
\sk
For a bounded complex of mixed Hodge modules $\MM$, we can define the {\it Hirzebruch class}
$$T_{y*}(\MM)\in\HH_{\ssb}(X)[y,y^{-1}],$$
by using the filtered de Rham complex and the Todd class transformation \cite{BFM}, see (1.1) below. This can be lifted naturally to the {\it spectral Hirzebruch class}
$$T^{\,\spr}_{\y*}(\MM,T_s)\in\HH_{\ssb}(X)\bl[\y^{\,1/e},\y^{\,-1/e}\br],$$
if $\MM$ has an action $T_s$ of finite order $e$, where $\y=-y$, see (1.3) below (and also \cite[Remark 1.3(4)]{CMSS}). This is defined by extending the definition of the ``dual" ${\rm Sp}'(f,x)$ of the Steenbrink spectrum ${\rm Sp}(f,x)$ in \cite[Section 2.1]{ste}. Note that the former is called the Hodge spectrum in the definition before \cite[Corollary 6.24]{DL}, see also \cite[Section 6.1]{GLM}.
(In the case $X$ is a point, the spectral Hirzebruch class is identified with the Hodge spectrum as is explained in \cite[Remark 3.7]{CMSS}.)
We need this refinement of Hirzebruch classes, since there is a {\it shift} of the Hodge filtration $F$ in the Thom-Sebastiani theorem for filtered $\D$-modules depending on the eigenvalues of the Milnor monodromy (see Theorem~(3.2) below).
\sk
By definition we have the relation
$$\bl[T^{\,\spr}_{\y*}(\MM,T_s)\br]^{\rm int}=T_{y*}(\MM)\q\h{in}\,\,\,\,\HH_{\ssb}(X)[y,y^{-1}],
\leqno(1)$$
where $[*]^{\rm int}$ is defined to be the tensor product of the $\Q$-linear morphism
$$\Q\bl[\y^{\,1/e},\y^{\,-1/e}\br]\ni\msum_{i\in\Z}\,a_i\,\y^{\,i/e}\mapsto\msum_{i\in\Z}\,a_i(-y)^{[i/e]}\in\Q[y,y^{-1}],$$
with $[i/e]$ the integer part of $i/e$, and $a_i\in\Q$ ($i\in\Z$).
(This corresponds to forgetting the action of $T_s$.)
\sk
Let $s'_1,\dots,s'_n$ be sufficiently general sections of $L$ with $n:=\dim Y$, where $Y$, $X$, $L$ are as in the beginning of the introduction. (We have to choose a generic hyperplane section in order to get a one-parameter deformation of a hypersurface. For an inductive argument, we have to choose many.) Put
$$\Si_0:=\Si_X:=\Sing\,X.$$
For $j\in[1,n]$, set
$$X'_j:=s^{\prime\,-1}_j(0),\q U_j:=Y\setminus X'_j,\q \Si_j:=\mcap_{k=1}^j\,X'_k\cap\Si_X.$$
Here the assumption that the $s'_j$ are ``sufficiently general" means that $X'_j$ is transversal to any stratum of a Whitney stratification of $\Si_{j-1}$ by increasing induction on $j$. Put
$$r:=\max\bl\{j\mid\Si_j\ne\emptyset\br\}\,(\les\dim\Si_X\les\dim X=n-1).$$
Note that $r=\dim\Si_X$ if $m\ges 2$, and $r=0$ if $m=1$. For $j\in[0,r]$, set
$$Z^{\circ}_j:=\mcap_{k=1}^j\,X'_k\cap U_{j+1},\,\,\,\Si^{\circ}_j:=\Si_j\setminus X'_{j+1}\,\bl(=Z^{\circ}_j\cap\Si_X\br),\,\,\,f'_j:=(s/s_{j+1}^{\prime\,m})|_{Z^{\circ}_j}.$$
\sk
In this paper we show the following.
\msn
{\bf Theorem~1.} {\it There is a localized class $M_y(X)\in\HH_{\ssb}(\Si_X)[y]$, which is called the Hirzebruch-Milnor class of $X$, and satisfies
$$T^{\,\vir}_{y*}(X)-T_{y*}(X)=(i_{\Si_X,X})_*M_y(X)\q\q\h{with}$$
\vskip-5mm
$$M_y(X)=\msum_{j=0}^r\,\bl[\,T^{\,\spr}_{\y*}\bl((i_{\Si^{\circ}_j,\Si_X})_{!\,}\varphi_{f'_j}\Q_{h,Z^{\circ}_j},T_s{\br)}\cdot\bl(-\msum_{i=1}^{m-1}\,\y^{\,i/m}\br){}^j\,\br]^{\rm int},
\leqno(2)$$
\vskip-5mm
$$T^{\,\spr}_{\y*}\bl((i_{\Si^{\circ}_j,\Si_j})_{!\,}\varphi_{f'_j}\Q_{h,Z^{\circ}_j},T_s{\br)}\in\HH_{\ssb}(\Si_j)\bl[\y^{\,1/e}\br]\q(j\in[0,r]),
\leqno(3)$$
where $i_{A,B}:A\into B$ denotes the inclusion for $A\subset B$ in general, $T_s$ in $(2)$, $(3)$ is the semisimple part of the monodromy, and $e$ is replaced by a multiple of it if necessary.}
\ms
In the case $\dim\Si_X=0$, we have $r=0$ and $\Si^{\circ}_0=\Si_X$. So the formula is quite simple using only the {\it Hirzebruch-Milnor classes} of the isolated singularities of $X$ by (1). In the case $\dim\Si_X=1$, we have $r=1$, and the right-hand side of (2) consists of two terms. The first term is the {\it Hirzebruch class} of the zero extension of the vanishing cycle Hodge module over the first generic hyperplane section by using (1). The second term is given by taking the {\it integer part} of the {\it spectral Hirzebruch-Milnor classes} of the isolated singularities of the first generic hyperplane section of $X$ which are {\it multiplied} by the fractional polynomial $-\msum_{i=1}^{m-1}\,\y^{\,i/m}$. (The situation is similar for the case $\dim\Si_X\ges 2$.)
\sk
Let $\ab=(a_1,\dots,a_n)$ be a set of sufficiently general non-zero complex numbers with $|a_j|$ sufficiently small. For $j\in[0,n-1]$, put
$$s_{\ab,j}:=s-\msum_{k=1}^j\,a_ks^{\prime\,m}_k,\q X_{\ab,j}:=s_{\ab,j}^{-1}(0),\q f_{\ab,j}:=\bl(s_{\ab,j}/s_{j+1}^{\prime\,m})|_{U_{j+1}}.$$
Note that $\Si_j=\Sing\,X_{\ab,j}$ if $m\ges 2$ (since the $s'_j$, $a_j$ are sufficiently general).
\sk
Theorem~1 is a consequence of the following two theorems.
\msn
{\bf Theorem~2.} {\it We have the Hirzebruch-Milnor class $M_y(X)\in\HH_{\ssb}(\Si_X)[y]$, satisfying}
$$T^{\,\vir}_{y*}(X)-T_{y*}(X)=(i_{\Si_X,X})_*M_y(X)\q\q\h{\it with}$$
\vskip-7mm
$$M_y(X)=\msum_{j=0}^r\,T_{y*}\bl((i_{\Si^{\circ}_j,\Si_X})_{!\,}\varphi_{f_{\ab,j}}\Q_{h,U_j}\br),
\leqno(4)$$
\msn
{\bf Theorem~3.} {\it For $j\in[0,r]$, we have
$$T^{\,\spr}_{\y*}\bl((i_{\Si^{\circ}_j,\Si_j})_{!\,}\varphi_{f'_j}\Q_{h,Z^{\circ}_j},T_s{\br)}\in\HH_{\ssb}(\Si_j)\bl[\y^{\,1/e}\br],
\leqno(5)$$
and there are equalities in $\HH_{\ssb}(\Si_j)\bl[\y^{\,1/e}\br]:$
$$\aligned&T^{\,\spr}_{\y*}\bl((i_{\Si^{\circ}_j,\Si_j})_{!\,}\varphi_{f_{\ab,j}}\Q_{h,U_j},T_s{\br)}\\&=T^{\,\spr}_{\y*}\bl((i_{\Si^{\circ}_j,\Si_j})_{!\,}\varphi_{f'_j}\Q_{h,Z^{\circ}_j},T_s{\br)}\cdot\bl(-\msum_{i=1}^{m-1}\,\y^{\,i/m}\br)^j,\endaligned
\leqno(6)$$
where $e$ is replaced by a multiple of it if necessary.}
\ms
The assertion (4) may be viewed as an inductive formula, since we have the following.
\msn
{\bf Proposition~1.} {\it For $j\in[0,r]$, there are equalities in $\HH_{\ssb}(X_{\ab,j})[y]:$}
$$\lim_{a_{j+1}\to 0}T_{y*}(X_{\ab,j+1})-T_{y*}(X_{\ab,j})=T_{y*}\bl((i_{\Si^{\circ}_j,X_{\ab,j}})_{!\,}\varphi_{f_{\ab,j}}\Q_{h,U_j}\br).
\leqno(7)$$
\ms
The limit in (7) is defined by using the nearby cycle functor $\psi$ for mixed Hodge modules, see (2.3.1) below.
Note that $X_{\ab,0}=X$, and $X_{\ab,r+1}$ is smooth. (If $r=n-1=\dim X$, then $s_{\ab,n}$ and $X_{\ab,n}$ can be defined in the same way as above by choosing a sufficiently general $s'_{n+1}$, and $X_{\ab,n}$ is smooth.) We assumed $m=1$ in \cite{MSS1}, where the formula was rather simple in the hypersurface case (since $r=0$ and $X_{\ab,1}$ is smooth, if $m=1$).
\sk
For the proof of Theorem~3, we need the Thom-Sebastiani theorem for the underlying filtered $\D$-modules of vanishing cycles (see Theorem~(3.2) below). Its relatively simple proof (using the algebraic microlocalization as is mentioned in \cite[Remark 4.5]{mic}) is given in \cite{MSS3}.
Note that, in the proof of Theorem~3, we can apply the Thom-Sebastiani theorem only by restricting to $U_j=Y\setminus X'_j$ and moreover only after passing to the normal bundle of $Z^{\circ}_j\subset U_{j+1}$ by using the deformation to the normal bundle, see Section~3 below.
\sk
The above Thom-Sebastiani theorem also implies the following Thom-Sebastiani type theorem for the {\it localized spectral Hirzebruch-Milnor classes} of $X:$
$$M_{\y}^{\,\spr}(X):=T_{\y}^{\,\spr}(\varphi_f\Q_{h,Y},T_s)\in\HH_{\ssb}(\Si_X)\bl[\y^{\,1/e}\br],$$
in the case where $X=f^{-1}(0)$ with $f$ a non-constant function on a smooth complex variety or a connected complex manifold $Y$ (that is, $f\in\Gamma(Y,\OO_Y)\setminus\C$), and $\Si_X:=\Sing\,X$.
\msn
{\bf Theorem~4.} {\it Let $X_a:=f_a^{-1}(0)$ with $f_a$ a non-constant function on a smooth complex variety or a connected complex manifold $Y_a$ {$(a=1,2)$}. Set $X:=f^{-1}(0)\subset Y:=Y_1\times Y_2$ with $f:=f_1+f_2$. Put $\Si_a:=\Sing\,X_a$ {$(a=1,2)$}. Then we have the equality
$$M_{\y}^{\,\spr}(X)=-M_{\y}^{\,\spr}(X_1)\boxtimes M_{\y}^{\,\spr}(X_2)\q\h{in}\,\,\,\,\HH_{\ssb}(\Si_X)\bl[\y^{\,1/e}\br],
\leqno(8)$$
by replacing {$($if necessary$)$} $Y_a$ with an open neighborhood of $X_a$ {$(a=1,2)$} so that $\Si_X=\Si_{X_1}\times\Si_{X_2}$, where $\boxtimes$ is defined by using cross products or K\"unneth maps.}
\ms
This assertion holds at the level of Grothendieck groups (more precisely, in $K_0(\Si_X)\bl[\y^{\,1/e}\br]$, see (3.5.1) below). Here $\Si_X\ne\Si_{X_1}\times\Si_{X_2}$ if and only if there are non-zero critical values $c_a$ of $f_a$ ($a=1,2$) with $c_1+c_2=0$. Note that $X$ is always non-compact even if $X_1,X_2$ are compact, see Remarks~(3.5)(iii) below.
In the case of isolated hypersurface singularities, Theorem~4 is essentially equivalent to the Thom-Sebastiani theorem for the spectrum as in \cite{ScSt}, \cite{Va}.
We can calculate the localized spectral Hirzebruch-Milnor class $M_{\y}^{\,\spr}(X)$ by improving the arguments in \cite[Section 5]{MSS1} (keeping track of the action of $T_s$), see also \cite[Remark 1.3]{CMSS}
\sk
It is known (see \cite[Theorem 3.2]{CMSS}, \cite[Corollary 3.12]{Sch5}) that we have
$$(i_{\Si_X,X})_*\bl[M_{\y}^{\,\spr}(X)\br]^{\rm int}{}={}(i_{\Si_X,X})_*M_y(X){}=T_{y*}^{\,\vir}(X)-T_{y*}(X).
\leqno(9)$$
Here the first equality follows from the definition, and the last one from the short exact sequence associated with the nearby and vanishing cycles together with \cite{Sch4} (or \cite[Proposition 3.3]{MSS1}) and \cite[Theorem 7.1]{Ve}. (We may also need \cite[Proposition 1.3.1]{MSS1} to show some compatibility of definitions, see the remark about $T_{y*}^{\,\vir}(X)$ after (1.1.9) below.)
\sk
Specializing to $y=-1$ (that is, $\y=1$), Theorems~1, 2 and 3 imply the corresponding assertions for the Chern classes. Indeed, $T_{y*}(X)$ and $T^{\,\vir}_{y*}(X)$ respectively specialize at $y=-1$ to the MacPherson-Chern class $c(X)$ (see \cite{Mac}) and the virtual Chern class $c^{\vir}(X)$ (called the Fulton or Fulton-Johnson class, see \cite{Fu}, \cite{FJ}) with rational coefficients, see \cite[Proposition 5.21]{Sch3}. The specialization of Theorem~4 to $y=-1$ is already known, see \cite[Section 4]{OY}, where a different sign convention is used.
\sk
We can moreover deduce a rather simpler formula for {\it localized Milnor classes} as follows: In the notation of Theorem~1, we define the {\it reduced Euler-Milnor function} $\EM_X$ on $\Si_X$ by
$$\EM_X(x)=\chi(F_{h_x})-1\q(x\in\Si_X),$$
where $F_{h_x}$ is a Milnor fiber of a local defining function $h_x$ of $X$ at $x$. (This cannot be generalized to spectrum functions, since the information of local system monodromies is also needed in this case, see \cite{MSS2}.) From Theorem~1 we can deduce the following by specializing to $y=-1$ and using \cite[Proposition 5.21]{Sch3}, where we assume $m\ges 2$.
\msn
{\bf Corollary~1.} {\it In the above notation and assumption, we have the localized Milnor class $M(X)\in\HH_{\ssb}(\Si_X)$ satisfying $c^{\vir}(X)-c(X)=(i_{\Si_X,X})_*M(X)$, and}
$$\aligned M(X)&=c(\EM_X)-\msum_{j=1}^r\,m(1-m)^{j-1}(i_{\Si_j,\Si_X})_*c(\EM_X|_{\Si_j})\\&=c(\EM_X)\cap\bl(1-\msum_{j=1}^r\,m(1-m)^{j-1}\bl(\tfrac{c_1(L)}{1+c_1(L)}\br)^j\,\br)\,\in\HH_{\ssb}(\Si_X).\endaligned
\leqno(10)$$
\ms
Here $c$ denotes MacPherson's Chern class transformation \cite{Mac}. For the proof of the second equality of (10), we use an analogue of Verdier's Riemann-Roch theorem for the MacPherson Chern classes under transversal pullbacks, see \cite[Corollary 0.1]{Sch1} or \cite[Corollary 2.7]{Sch6}. Restricting to the degree 0 part, Corollary~1 implies the following assertion about Euler numbers.
\msn
{\bf Corollary~2.} {\it In the above notation, we have
$$\chi(X')-\chi(X)=(a_{\Si_X})_*(\EM_X)\,-\,\msum_{j=1}^r\,m(1-m)^{j-1}(a_{\Si_j})_*(\EM_X|_{\Si_j}),
\leqno(11)$$
with $X'$ a smooth hypersurface of degree $m$ in $Y\,\,($with respect to the ample line bundle $L)$.}
\ms
Here $(a_{\Si_X})_*(\EM_X):=\sum_S\chi(S)\,c_S$ with $\{S\}$ a (smooth) stratification of $\Si_X$ satisfying $\EM_X|_{S}=c_S\in\Z$ (see \cite{Mac}), and similarly for $(a_{\Si_j})_*(\EM_X|_{\Si_j})$. Corollaries~1 and 2 are well-known in the isolated singularity case (that is, $r=0$) where the right-hand side of (10) and (11) is given by the sum of Milnor numbers. (The reader may assume $Y=\PP^n$ and $L=\OO_{\PP^n}(1)$ if he likes.) There is a similar formula
$$\aligned M(X)&=c(\EM_X)-(i_{\Si'_X,\Si_X})_*c(\EM_X|_{\Si'_X})\\&=c(\EM_X)(1+m\,c_1(L))^{-1}\in\HH_{\ssb}(\Si_X),\endaligned
\leqno(12)$$
where $\Si'_X:=\Si_X\cap X'$ with $X'$ sufficiently general. The first equality of (12) follows from \cite[Theorem 5.3]{PP}, and the second is shown by an argument similar to the proof of the second equality of (10). The obtained equality is closely related with \cite[Theorem 0.2]{PP} and \cite[Formula 14]{Sch1}. (We can show directly the equality between the last terms of (10) and (12), see Remark~(2.7) below.) Taking the degree 0 part of (12), we also get
$$\chi(X')-\chi(X)=(a_{\Si_X})_*(\EM_X)-(a_{\Si'_X})_*(\EM_X|_{\Si'_X}).
\leqno(13)$$
In this case, however, we have to calculate the intersection of $\Si_X=\Sing\,X$ with a sufficiently general hypersurface $X'$ of degree $m$, and this can be rather complicated if $m$ is quite large. Sometimes the calculation may be simpler if it is enough to consider only the restrictions to intersections of sufficiently general {\it hyperplanes}, see for instance Example~(2.8) below.
\sk
As for an application to du Bois singularities, we have the following.
\msn
{\bf Theorem~5.} {\it Let $X$, $Y$ be as in Theorem~$1$ with $X$ reduced, or $X$ be a reduced hypersurface in a smooth complex algebraic variety or a complex manifold $Y$ defined by $f\in\Gamma(Y,\OO_Y)$. If $X$ has only Du Bois singularities, then
$$M_0(X):=M_y(X)|_{y=0}=0\q\h{in}\,\,\,\,\HH_{\ssb}(\Si_X).
\leqno(14)$$
The converse holds if $\Si_X={\rm Sing}\,X$ is a projective variety.
More precisely, the converse holds if $(14)$ holds for $(i_{\Si_X,\PP^N})_*M_0(X)$ in $\HH_{\ssb}(\PP^N)$, where $\PP^N$ is a projective space containing $\Si_X$.}
\ms
The first assertion of Theorem~5 is already known if $\HH_{\ssb}(\Si_X)$ in (14) is replaced with $\HH_{\ssb}(X)$, and $M_0$ with $(i_{\Si_X,X})_*M_0$ at least in the second case where $X$ is defined by a global function $f$, since we have
$$T^{\,\vir}_{y*}(X)|_{y=0}=td_*(\OO_X),$$
see \cite[p.~6]{BSY}, \cite[p.~2619]{CMSS}, \cite[Corollary 2.3]{Sch5}. Its converse also holds if $X$ is projective.
In the isolated singularity case, Theorem~5 is related to \cite[Theorem 3.12]{St3}, \cite[Theorem 6.3]{Is} in the case where the smoothing is a base change of a smoothing with total space nonsingular, see (4.8) below.
\sk
Spectral Hirzebruch-Milnor classes are related to jumping coefficients as follows.
\msn
{\bf Proposition~2.} {\it Let $X$ be a hypersurface in a smooth complex algebraic variety {\rm(}or a complex manifold{\,$)$} $Y$ defined by $f\in\Gamma(Y,\OO_Y)$. Let $\al\in(0,1)$. Then $\al\in\JC(f)$ if
$$M^{\rm sp}_{\y}(X)|_{\y^{\al}}\ne 0\q\h{in}\,\,\,\,\HH_{\ssb}(\Si_X),
\leqno(15)$$
where $M^{\rm sp}_{\y}(X)|_{\y^{\al}}$ is the coefficient of $\y^{\al}$ in $M^{\rm sp}_{\y}(X)\in\HH_{\ssb}(\Si_X)\bl[\y^{1/e}\br]$.
The converse holds if $\Si_X={\rm Sing}\,X$ is a projective variety.
More precisely, if $\al\in\JC(f)\cap(0,1)$, then $(15)$ holds for the image of $M^{\rm sp}_{\y}(X)|_{\y^{\al}}$ in $\HH_{\ssb}(\PP^N)$, where $\PP^N$ is a projective space containing $\Si_X$.}
\ms
Here $\JC(f)$ is the set of {\it jumping coefficients} consisting of $\al\in\Q$ with $\G(\al X)\ne 0$, and the $\G(\al X)$ are graded quotients of multiplier ideal sheaves $\J(\al X)$, see \cite{La} (and also \cite[Section 2]{MSS3}). Proposition~2 is an immediate consequence of the duality for nearby cycle functors \cite{dual} combined with \cite[Theorem 0.1]{BS1} (see also (4.3.10) below).
In the isolated singularity case, it is closely related to \cite[Corollary on p.~258]{Bu}. We have by the second equality of (4.5.2) below
$$M_0(X)=\mopl_{\al\in(0,1)}\,M^{\rm sp}_{\y}(X)|_{\y^{\al}}\q\h{in}\,\,\,\,\HH_{\ssb}(\Si_X).
\leqno(16)$$
\sk
Let $\lct(f)$ be the {\it log canonical threshold} of $f$, which is by definition the minimal jumping coefficient. Then Proposition~2 may be viewed as a refinement of Theorem~5 modulo the assertion that we have $\lct(f)=1$ if and only if $X=f^{-1}(0)$ is reduced, and has only Du Bois singularities, see \cite[Theorem 0.5]{fil}, \cite[Corollary 6.6]{KoSc}, and (4.3.9) below.
\sk
The first named author is partially supported by NSF and NSA. The second named author is partially supported by Kakenhi 15K04816. The third named author is supported by the SFB 878 ``groups, geometry and actions''.
\sk
In Section~1 we review some basics of Hirzebruch classes, introduce spectral Hirzebruch-Milnor classes showing some of their properties, and explain the topological filtration on the Grothendieck groups. In Section~2 we prove Theorem~2 and Proposition~1 by using the short exact sequences associated with nearby and vanishing cycle functors of mixed Hodge modules. In Section~3 we prove Theorems~3 and 4 after explaining the Thom-Sebastiani theorem for underlying filtered $\D$-modules of constant mixed Hodge modules. In Section~4 some relations with rational and Du Bois singularities are explained.
\bs\bs
\vbox{\centerline{\bf 1. Hirzebruch characteristic classes}
\bsn
In this section we review some basics of Hirzebruch classes, introduce spectral Hirzebruch-Milnor classes showing some of their properties, and explain the topological filtration on the Grothendieck groups.}
\msn
{\bf 1.1.~Hirzebruch classes.} For a complex algebraic variety $X$, we set in this paper
$$\HH_k(X):=H^{\BM}_{2k}(X,\Q)\q\h{or}\q\CH_k(X)_{\Q}.$$
Let $\MHM(X)$ be the abelian category of mixed Hodge modules on $X$ (see \cite{mhp}, \cite{mhm}). For $\MM\in D^b\MHM(X)$, we have the homology Hirzebruch characteristic class defined by
$$T_{y*}(\MM):=td_{(1+y)*}\bl(\DR_y[\MM]{\br)}\in\HH_{\ssb}(X)[y,y^{-1}].
\leqno(1.1.1)$$
Here, setting $F^p=F_{-p}$ so that $\Gr_F^p=\Gr^F_{-p}$, we have
$$\DR_y[\MM]:=\msum_{i,p}\,(-1)^i\,\bl[\Hc^i\Gr_F^p\DR(\MM)\br]\,(-y)^p\in K_0(X)[y,y^{-1}],
\leqno(1.1.2)$$
with
$$td_{(1+y)*}:K_0(X)[y,y^{-1}]\to\HH_{\ssb}(X)\bl[y,\h{$\frac{1}{y(y+1)}$}\br]
\leqno(1.1.3)$$
defined by the composition of the scalar extension of the Todd class transformation
$$td_*:K_0(X)\to\HH_{\ssb}(X),$$
(see \cite{BFM}, where $td_*$ is denoted by $\tau$) with the multiplication by $(1+y)^{-k}$ on $\HH_k(X)$ (see \cite{BSY}).
Note that $\Hc^i\Gr_F^p\DR(\MM)=0$ for $|p|\gg 0$ (see for instance \cite[2.2.10.5]{mhp}).
We have the last inclusion in (1.1.1), that is, $T_{y*}(\MM)\in\HH_{\ssb}(X)[y,y^{-1}]$, by \cite[Proposition 5.21]{Sch3}.
\sk
In this paper, we denote $a_X^*\Q_h\in D^b\MHM(X)$ by $\Q_{h,X}$ as in \cite{MSS1}, where $\Q_h$ denotes the trivial $\Q$-Hodge structure of rank $1$ and type $(0,0)$ (see \cite{De1}), and $a_X:X\to pt$ is the structure morphism, see \cite{mhm}.
The {\it homology Hirzebruch characteristic class} of $X$ is defined by
$$T_{y*}(X):=T_{y*}(\Q_{h,X})\in\HH_{\ssb}(X)[y].
\leqno(1.1.4)$$
The last inclusion is reduced to the $X$ smooth case by using a stratification of $X$ together with smooth partial compactifications of strata over $X$. Then it follows from the relation with the cohomology Hirzebruch class, see \cite{BSY}.
\sk
In the $X$ smooth case, we have
$$\DR_y[X]=\La_y[T^*X],
\leqno(1.1.5)$$
with
$$\La_y[V]:=\msum_{p\ges 0}\,[\La^pV]\,y^p\q\h{for a vector bundle}\,\,\,V.
\leqno(1.1.6)$$
\sk
In the case $X$ is a complete intersection in a smooth complex algebraic variety $Y$, the {\it virtual Hirzebruch characteristic class} $T^{\,\vir}_{y*}(X)$ can be defined like the virtual genus in \cite{Hi} (see \cite[Section 1.4]{MSS2}) by
$$T^{\,\vir}_{y*}(X):=td_{(1+y)*}\DR^{\vir}_y[X]\in\HH_{\ssb}(X)[y],
\leqno(1.1.7)$$
where $\DR^{\vir}_y[X]$ is the image in $K_0(X)[[y]]$ of
$$\La_y(T^*_{\vir}X)=\La_y[T^*Y|_X]/\La_y[N^*_{X/Y}]\in K^0(X)[[y]],
\leqno(1.1.8)$$
and belongs to $K_0(X)[y]$ (see for instance \cite[Proposition 3.4]{MSS1}).
Here $K^0(X)$, $K_0(X)$ are the Grothendieck group of locally free sheaves of finite length and that of coherent sheaves respectively.
We denote by $T^*Y$ and $N^*_{X/Y}$ the cotangent and conormal bundles respectively. The virtual cotangent bundle is defined by
$$T^*_{\vir}X:=[T^*Y|_X]-[N^*_{X/Y}]\in K^0(X).
\leqno(1.1.9)$$
Here $N^*_{X/Y}$ in the non-reduced case is defined by the locally free sheaf $\I_X/\I_X^2$ on $X$ with $\I_X\subset\OO_Y$ the ideal of $X\subset Y$.
\sk
Note that the above definition of $T^{\,\vir}_{y*}(X)$ is compatible with the one in \cite{CMSS} by \cite[Proposition 1.3.1]{MSS1}.
\msn
{\bf Remarks~1.2.} (i) Let $(M,F)$ be a filtered {\it left} $\D$-module on a smooth variety $X$ of dimension $d_X$. The filtration $F$ of the de Rham complex $\DR_X(M,F)$ is defined by
$$\DR_X(M,F)^i=\Omega_X^{i+d_X}{\otimes_{\OO_X}}(M,F[-i-d_X])\q\q\bl(i\in[-d_X,0]\br),
\leqno(1.2.1)$$
where $\DR_X(M,F)^i$ denotes the $i\,$th component of $\DR_X(M,F)$. Recall that $F_p=F^{-p}$ and $F[m]_p=F_{p-m}$ for $p,m\in\Z$.
\sk
This is compatible with the definition of $\DR_X$ for filtered {\it right} $\D_X$-modules as in \cite{mhp}; that is, we have the canonical isomorphism
$$\DR_X(M,F)=\DR_X\bl((\Omega_X^{d_X},F){\otimes_{\OO_X}}(M,F)\br),
\leqno(1.2.2)$$
where $(\Omega_X^{d_X},F){\otimes_{\OO_X}}(M,F)$ is the filtered right $\D_X$-module corresponding to a filtered $\D_X$-module $(M,F)$, and the filtration $F$ on $\Omega_X^{d_X}$ is shifted by $-d_X$ so that
$$\Gr^F_p\Omega_X^{d_X}=0\q(\,p\ne-d_X),
\leqno(1.2.3)$$
\sk
Under the direct image functor $i_*^{\D}$ for filtered $\D$-modules with $i:X\into Y$ a closed embedding of smooth varieties, the filtration $F$ of a filtered {\it left} $\D$-module $(M,F)$ is shifted by the codimension $r:=\dim Y-\dim X$; more precisely
$$i_*^{\D}M=i_*M[\dd_1,\dots,\dd_r]\q\h{with}\q F_p(i_*^{\D}M)=\msum_{\nu\in\N^r}\,i_*\bl(F_{p-|\nu|-r}M\otimes\dd^{\nu}\br),
\leqno(1.2.4)$$
where $\dd^{\nu}:=\prod_{j=1}^r\dd_j^{\nu_j}$ with $\dd_j:=\dd/\dd y_j$ for local coordinates $y_i$ of $Y$ with $X=\mcup_{i\les r}\{y_i=0\}$.
This shift comes from (1.2.3) since there is no shift of filtration for filtered right $\D$-modules.
(Globally we have to twist the right-hand side of (1.2.4) by $\om_{X/Y}$.)
Because of this shift of the filtration $F$, we have the canonical isomorphisms of $\OO_Y$-modules
$$i_*\Hc^j\Gr^F_p\DR_X(M)=\Hc^j\Gr^F_p\DR_Y(i_*^{\D}M).
\leqno(1.2.5)$$
\sk
If $(M,F)$ is the underlying filtered left $\D_X$-module of a mixed Hodge module $\M$, then there is an equality in $K_0(X)[y,y^{-1}]$:
$$\aligned&\msum_{i,p}\,(-1)^i\,\bl[\Hc^i\Gr_F^p\DR(M)\br](-y)^p\\&=\msum_{i,p}\,(-1)^i\,\bl[\Omega_X^{i+d_X}{\otimes_{\OO_X}}\Gr_F^pM\br](-y)^{p+i+d_X}.\endaligned
\leqno(1.2.6)$$
This is related to the right-hand side of (1.1.2), and follows directly from (1.2.1).
\ms
(ii) Let $Z$ be a locally closed smooth subvariety of a smooth variety $X$ with $i_Z:Z\into X$ the canonical inclusion. For $\MM\in D^b\MHM(X)$, set
$$\M_Z^j:=H^ji_Z^*\MM\in\MHM(Z)\q(j\in\Z).$$
(Here $H^j:D^b\MHM(Z)\to\MHM(Z)$ is the canonical cohomology functor.)
Assume this mixed Hodge module is a variation of mixed Hodge structure $\HH_Z^j$ on $Z$ for any $j\in\Z$; more precisely, its underlying $F$-filtered left $\D_Z$-module is a locally free $F$-filtered $\OO_Z$-module which underlies $\HH_Z^j$. (Note that there is a shift of the weight filtration $W$ by $d_Z$.) For $z\in Z$, we denote by $i_{z,Z}:\{z\}\into Z$ and $i_z:\{z\}\into X$ the canonical inclusions. We have $\M_z^j$, $\HH_z^j$ by applying the above argument to the inclusion $i_z:\{z\}\into X$. Here $\M_z^j$ is naturally identified with $\HH_z^j$ (including the weight filtration $W$) since $z$ is a point. The relations between these are given by
$$H^{-d_Z}i_{z,Z}^*\M_Z^j=i_{z,Z}^*\HH_Z^j=\M_z^{j-d_Z}=\HH_z^{j-d_Z}.
\leqno(1.2.7)$$
Here $H^ki_{z,Z}^*\M_Z^j=0$ ($k\ne -d_Z$), and $i_{z,Z}^*\HH_Z^j$ is the restriction to $z$ as a variation of mixed Hodge structure.
Note that there is {\it no shift} of the filtration $F$ in (1.2.7).
(Consider, for instance, the case $\MM$ is the constant mixed Hodge module $\Q_{h,X}[d_X]$, where $\M_Z^{-r}=\Q_{h,Z}[d_Z]$ with $r={\rm codim}_XZ$ and $\M_z^{-d_X}=\Q_{h}$.)
This follows from the definition of the nearby cycle functors $\psi_{z_i}$ in \cite{mhp}, \cite{mhm}, where the $z_i$ are local coordinates of $Z$. The functor $i_{z,Z}^*$ can be given in this case by the iteration of the mapping cones of the canonical morphisms
$$\can:\psi_{z_i,1}\to\varphi_{z_i,1},$$
and the vanishing cycle functors $\varphi_{z_i,1}$ vanish for smooth mixed Hodge modules on $Z$. Here smooth means that their underlying $\Q$-complexes are local systems on $Z$ shifted by $d_Z$. (Note that the last property implies the shift of indices in (1.2.7).)
\msn
{\bf 1.3.~Spectral Hirzebruch classes.} We denote by $\MHM(X,T_s)$ the abelian category of mixed Hodge modules $\M$ on a smooth variety $X$ (or more generally, on a variety embeddable into a smooth variety $X$) such that $\M$ is endowed with an action of $T_s$ of finite order. (For instance, $\M=\varphi_{f_{\ab,j}}\Q_{h,U_j}$ with $T_s$ the semisimple part of the monodromy in the notation of the introduction.)
For $(\M,T_s)\in\MHM(X,T_s)$, let $(M,F)$ be the underlying filtered left $\D_X$-module. Since $T_s$ has a finite oder $e$, we have the canonical decomposition
$$(M,F)=\msum_{\la\in\mu_e}\,(M_{\la},F),
\leqno(1.3.1)$$
such that $T_s=\la$ on $M_{\la}\subset M$, where $\mu_e:=\{\la\in\C\mid\la^e=1\}$.
We define the {\it spectral Hirzebruch class} by
$$T^{\,\spr}_{\y*}(\M,T_s):=td_{(1-\y)*}\bl(\DR_{\y}[\M,T_s]{\br)}\in\HH_{\ssb}(X)\bl[\y^{\,1/e},\h{$\frac{1}{\y(\y-1)}$}\br],
\leqno(1.3.2)$$
with
$$\aligned\DR_{\y}[\M,T_s]:=\msum_{i,p,\la}\,(-1)^i\,\bl[\Hc^i\Gr_F^p\DR(M_{\la})\br]\,\y^{\,p+\ell(\la)}\q\q\q&\\ \h{in}\,\,\,\,K_0(X)[\y^{\,1/e},\y^{\,-1/e}].&\endaligned
\leqno(1.3.3)$$
Here
$$\ell(\la)\in[0,1)\q\h{with}\q\exp(2\pi i\ell(\la))=\la,
\leqno(1.3.4)$$
and
$$td_{(1-\y)*}:K_0(X)[\y^{\,1/e},\y^{\,-1/e}]\to\HH_{\ssb}(X)\bl[\y^{\,1/e},\h{$\frac{1}{\y(\y-1)}$}\br]
\leqno(1.3.5)$$
is the scalar extension of the Todd class transformation $td_*:K_0(X)\to\HH_{\ssb}(X)$ followed by the multiplication by $(1-\y)^{-k}$ on $\HH_k(X)$ as in (1.1) (where $\y=-y$). Actually the class belongs to $\HH_{\ssb}(X)[\y^{\,1/e},\y^{\,-1/e}]$ by generalizing \cite[Proposition 5.21]{Sch3}, see Proposition~(1.4) below.
\sk
The above definition can be generalized to the $X$ singular case by using locally defined closed embeddings into smooth varieties, where the independence of locally defined closed embeddings follows from the isomorphism (1.2.5) (at the level of $\OO$-modules).
We can further generalize this definition to the case of $\MM\in D^b\MHM(X)$ endowed with an action of $T_s$ of finite order by applying the above arguments to each cohomology module $H^i\MM$ ($i\in\Z$).
\sk
The above arguments imply the transformation
$$T^{\,\spr}_{\y*}:K_0^{\rm mon}(\MHM(X))\to\mcup_{e\ges 1}\,\HH_{\ssb}(X)\bl[\y^{\,1/e},\y^{\,-1/e}\br],
\leqno(1.3.6)$$
which is functorial for proper morphisms by using the compatibility of $\DR$ (or rather $\DR^{-1}$) with the direct images by proper morphisms (see \cite[Section 2.3.7]{mhp}) together with the compatibility of $td_*$ with the pushforward by proper morphisms (see \cite{BFM}).
Here the left-hand side is the Grothendieck group of mixed Hodge modules on $X$ endowed with a finite order automorphism, see also \cite[Remark 1.3(4)]{CMSS}.
\msn
{\bf Proposition~1.4.} {\it In the notation of $(1.3.2)$, we have}
$$T^{\,\spr}_{\y*}(\M,T_s)\in\HH_{\ssb}(X)[\y^{\,1/e},\y^{\,-1/e}].
\leqno(1.4.1)$$
\msn
{\it Proof.} Let $\Sc$ be a stratification of $X$ such that $H^ij_S^*\M$ are variations of mixed Hodge structures on any strata $S\in\Sc$ for any $i\in\Z$, where $j_S:S\into X$ is the canonical inclusion.
By the same argument as in the proof of \cite[Proposition 5.1.2]{MSS1}, we have the equality
$$T^{\,\spr}_{\y*}(\M)=\msum_{S,i}\,(-1)^iT^{\,\spr}_{\y*}\bl((j_S)_!H^i(j_S)^*\M\br).
\leqno(1.4.2)$$
So the assertion is reduced to the case where $\M=j_!\M'$ with $\M'\in\MHM(X')$ a variation of mixed Hodge structure on an open subvariety $X'$ with $j:X'\into X$ the natural inclusion.
We may further assume that $D:=X\setminus X'$ is a divisor with simple normal crossings since the Todd class transformation $td_*$ and the de Rham functor $\DR$ commute with the pushforward or the direct image under a proper morphism (see \cite{BFM} for $td_*$).
\sk
Let $M^{>0}$ be the Deligne extension of the underlying $\OO_{X'}$-module of $\M'$ such that the eigenvalues of the residues of the logarithmic connections are contained in $(0,1]$. The action of $T_s$ is naturally extended to $M^{>0}$, and we have the canonical decomposition
$$M^{>0}=\mopl_{\la}\,M^{>0}_{\la}.$$
It follows from \cite[Proposition 3.11]{mhm} that each $M^{>0}_{\la}$ is identified with an $\OO_X$-submodule of $M_{\la}$, and there is a canonical filtered quasi-isomorphism
$$\DR_{X\langle D\rangle}(M^{>0}_{\la},F)\simto\DR_X(M_{\la},F),
\leqno(1.4.3)$$
where the left-hand side is the filtered logarithmic de Rham complex such that its $i\,$th component is given by
$$\DR_{X\langle D\rangle}(M^{>0}_{\la},F)^i=\Omega_X^{i+d_X}(\log D){\otimes_{\OO_X}}(M^{>0}_{\la},F[-i-d_X]).
\leqno(1.4.4)$$
As in (1.2.6), we get by (1.4.3) the equality in $K_0(X)[\y^{\,1/e},\y^{\,-1/e}]$:
$$\aligned&\q\msum_{p,i,\la}\,(-1)^i\,\bl[\Hc^i\Gr_F^p\DR_X(M_{\la})\br]\,\y^{\,p+\ell(\la)}\\&=\msum_{p,i,\la}\,(-1)^i\,\bl[\Omega_X^{i+d_X}(\log D){\otimes_{\OO_X}}\Gr_F^pM^{>0}_{\la}\br]\,\y^{\,p+i+d_X+\ell(\la)}.\endaligned
\leqno(1.4.5)$$
By \cite{BFM}, the Todd class transformation $td_*$ satisfies the following property:
$$td_*(V{\otimes}\,\xi)=ch(V)\cap td_*(\xi)\q\bl(V\in K^0(X),\,\,\xi\in K_0(X)\br).
\leqno(1.4.6)$$
By using the {\it twisted} Chern character
$$ch^{(1-\y)}:K^0(X)\ni V\mapsto\msum_{k\ges 0}\,ch^k(V)(1-\y)^k\in\HH^{\ssb}(X)[\y],
\leqno(1.4.7)$$
with $ch^k$ the $k\,$th component of the Chern character $ch$ (see \cite{Sch3}, \cite{Yo1}), (1.4.6) is extended to
$$\aligned td_{(1-\y)*}(V{\otimes}\,\widetilde{\xi})=ch^{(1-\y)}(V)\cap td_{(1-\y)*}(\widetilde{\xi})\\ \bl(V\in K^0(X),\,\,\widetilde{\xi}\in K_0(X)[\y^{\,1/e},\y^{\,-1/e}]\br).\endaligned
\leqno(1.4.8)$$
By applying this to the case $V=\Gr^F_pM^{>0}_{\la}$ in the right-hand side of (1.4.5), the assertion is now reduced to the case $\Gr^F_pM^{>0}_{\la}=\OO_X$, and to the case $\M=j_*\Q_{h,X'}[d_X]$, as in the proof of \cite[Proposition 5.21]{Sch3}.
\sk
The assertion is further reduced to the case $X'=X$ by using the weight filtration $W$ on the mixed Hodge module $j_*\Q_{h,X'}[d_X]$. Indeed, if $D_i$ ($i=1,\dots,m$) are the irreducible components of $D$, then
$$\Gr^W_k(j_*\Q_{h,X'}[d_X])=\mopl_{|I|=k-d_X}\Q_{h,D_I}(-|I|)[d_{D_I}]\q(I\subset\{1,\dots,m\}),
\leqno(1.4.9)$$
where $D_I:=\bigcap_{i\in I}D_i$, and $|I|={\rm codim}_XD_I$ so that $k=d_{D_I}+2|I|=d_X+|I|$.
\sk
In the case $\M=\Q_{h,X}[d_X]$, the assertion follows from the inclusion (1.1.4).
This finishes the proof of Proposition~(1.4).
\msn
{\bf Proposition~1.5.} {\it Let $D$ be a divisor on a smooth complex algebraic variety $X$. Set $U:=X\setminus D$ with $j:U\into X$ the natural inclusion. Let $f\in\Gamma(U,\OO_U)$. Then}
$$T^{\,\spr}_{\y*}(j_!\psi_f\Q_{h,U},T_s),\,\,T^{\,\spr}_{\y*}(j_!\varphi_f\Q_{h,U},T_s)\,\in\,\HH_{\ssb}(X)[\y^{\,1/e}].
\leqno(1.5.1)$$
\msn
{\it Proof.} By the argument in (1.4) together with (1.2.7) in Remark~(1.2)(ii) (where there is no shift of the filtration $F$), the assertion (1.5.1) can be reduced to the following:
$$\Gr_F^pH^j(F_{\!f,0},\C)=0\q\h{for}\,\,\, p<0.
\leqno(1.5.2)$$
Here $0$ can be any point of $f^{-1}(0)\subset U$, and we denote by $F_{\!f,0}$ the Milnor fiber of $f$ around $0$ so that
$$H^j(F_{\!f,0},\Q)=H^ji_0^*\psi_f\Q_{h,U}\q\q(j\in\Z),
\leqno(1.5.3)$$
with $i_0:\{0\}\into X$ the natural inclusion.
More precisely, (1.5.1) is reduced to (1.5.2) by using the inclusion (1.1.4) together with (1.4.8) and the non-negativity of the codimension $|I|$ in (1.4.9).
\sk
For the proof of (1.5.2), we use an embedded resolution $\pi:\X\to X$ of $f^{-1}(0)$, where we may assume that $D:=\pi^{-1}(0)\subset\X$ is a divisor (by taking a point-center blow-up first). The latter is a union of irreducible components of $\pi^{-1}f^{-1}(0)$, and is also a divisor with normal crossings. Set $\pi_0:=\pi|_D:D\to\{0\}$, and $\ft:=f\ssc\pi$. We have the canonical isomorphisms
$$H^ji_0^*\psi_f\Q_{h,X}=H^j(\pi_0)_*i_D^*\psi_{\ft}\Q_{h,\X}\q\q(j\in\Z),
\leqno(1.5.4)$$
since
$$\pi_*\psi_{\ft}\Q_{h,\X}=\psi_f\pi_*\Q_{h,\X}=\psi_f\Q_{h,X}\q\h{and}\q(\pi_0)_*\ssc i_D^*=i_0^*\ssc\pi_*.$$
It is well-known that the variation of mixed Hodge structures $\HH_Z^j$ in Remark~(1.2)(ii) (applied to $\MM=i_D^*\psi_{\ft}\Q_{h,\X}$ or $\psi_{\ft}\Q_{h,\X})$ are direct sums of locally constant variations of Hodge structures of type $(k,k)$ with $k\in\N$, where $Z$ is a stratum of the canonical stratification associated with the divisor with normal crossings $\pi^{-1}f^{-1}(0)$, see \cite{St1}, \cite{St2} (and also \cite[Proposition 3.5]{mhm}). So the assertion (1.5.2) follows. This finishes the proof of Proposition~(1.5).
\ms
We recall here the notion of topological filtration on the Grothendieck group (see \cite{SGA6}, \cite{Fu}) which will be used in Section~4.
\msn
{\bf 1.6.~Topological filtration.} Let $X$ be a complex algebraic variety, and $K_0(X)$ be the Grothendieck group of coherent sheaves on $X$. It has the topological filtration which is denoted in this paper by $G$ (in order to distinguish it from the Hodge filtration $F$), and such that $G_kK_0(X)$ is generated by the classes of coherent sheaves $\F$ with $\dim{\rm supp}\,\F\les k$, see \cite[Examples~1.6.5 and 15.1.5]{Fu}, \cite{SGA6}. It is known (see \cite[Corollary 18.3.2]{Fu}) that $td_*$ induces the isomorphisms
$$\aligned(td_*)_{\Q}&:K_0(X)_{\Q}\simto\mopl_k\,{\rm CH}_k(X)_{\Q},\\
\Gr^G_k(td_*)_{\Q}&:\Gr_k^GK_0(X)_{\Q}\simto{\rm CH}_k(X)_{\Q}.\endaligned
\leqno(1.6.1)$$
Moreover the inverse of the last isomorphism is given by $Z\mapsto [\OO_Z]$ for irreducible reduced closed subvarieties $Z\subset X$ with dimension $k$.
Here we set
$$\HH_{\ssb}(X):=\mopl_k\,{\rm CH}_k(X)_{\Q},$$
and use the filtration $G$ defined by
$$G_k\HH_{\ssb}(X):=\mopl_{j\les k}\,\HH_j(X).
\leqno(1.6.2)$$
This definition is valid also in the case $\HH_k(X):=H^{\rm BM}_{2k}(X,\Q)$. Here
$$\Gr^G_k(td_*)_{\Q}:\Gr_k^GK_0(X)_{\Q}={\rm CH}_k(X)_{\Q}\to H^{\rm BM}_{2k}(X,\Q)
\leqno(1.6.3)$$
is identified with the cycle class map.
\sk
If $X=\PP^N$, then $\HH_k(X)=\Q$ for $k\in[1,N]$. These are canonically generated by the classes of linear subspaces, and the $\Gr^G_k(td_*)_{\Q}$ are identified with the identity maps.
For any irreducible reduced closed subvariety $Z\subset\PP^N$ with $\dim Z=k$, we have the {\it positivity}\,:
$$\Gr^G_k(td_*)_{\Q}[Z]=\deg Z>0\q\h{in}\,\,\,\,\HH_k(\PP^N)=\Q.
\leqno(1.6.4)$$
Here $\deg Z$ is defined to be the intersection number of $Z$ with a sufficiently general linear subspace of the complementary dimension if $k>0$ (and $\deg Z=1$ if $k=0$).
Moreover, for any coherent sheaf $\F$ on $\PP^N$ with $\dim{\rm supp}\,\F=k$, we have also the {\it positivity}\,:
$$\Gr_k^G[\F]=\msum_{i=1}^r\,m_i\,\Gr_k^G[Z_i]>0\q\h{in}\,\,\,\,\Gr_k^GK_0(\PP^N)_{\Q}={\rm CH}_k(\PP^N)_{\Q}=\Q,
\leqno(1.6.5)$$
where the $Z_i$ are $k$-dimensional irreducible components of ${\rm supp}\F$, and $m_i\in\Z_{>0}$ are the multiplicity of $\F$ at the generic point of $Z_i$ ($i\in[1,r]$).
\msn
{\bf 1.7.~Hirzebruch class of a projective hyperplane arrangement.} Let $X\subset Y:=\PP^n$ be a projective hyperplane arrangement. Let $X_i$ ($i\in I_0$) be the irreducible components of $X$ (which are isomorphic to $\PP^{n-1}$). Set
$$\Sc(X):=\bl\{Z\subset X\,\,\big|\,\,Z=\mcap_{i\in I}\,X_i\,\,\,\h{for some}\,\,\,I\subset I_0\,\,\,\h{with}\,\,\,I\ne\emptyset\br\}$$
Put $\gamma(Z):={\rm codim}_YZ$. There are trivial local systems $L_Z$ of rank $r_Z$ on $Z\in\Sc(X)$ together with a quasi-isomorphism
$$\Q_X\simto K_X^{\ssb}\q\q\h{with}\q\q K_X^j:=\mopl_{\gamma(Z)=j+1}\,L_Z.
\leqno(1.7.1)$$
This follows for instance from \cite[Lemma 1.8]{BS2}. It implies that
$$\Gr_k^W(\Q_{h,X}[d_X])=\mopl_{\dim Z=k}\,L_{h,Z}[d_Z].
\leqno(1.7.2)$$
where $L_{h,Z}$ is a constant variation of Hodge structure of type $(0,0)$ on $Z$ such that the underlying $\Q$-local system is $L_Z$. (This is a generalization of the standard resolution of the constant sheaf on a variety with simple normal crossings by using the constant sheaves supported on intersections of irreducible components of the variety, see for instance \cite[Example 3.3]{St1}.)
\sk
The quasi-isomorphism (1.7.1) implies that $(-1)^{\gamma(Z)}r_Z$ coincides with the {\it M\"obius function} (see \cite[2.1]{OS}). (Specializing to $y=-1$, the above argument then seems to be closely related to a calculation of the Chern class of $\PP^n\setminus X$ in \cite{Al}.) By (1.7.2) the calculation of the Hirzebruch class of $X$ is reduced to that for projective spaces $Z\in\Sc(X)$, and the latter can be calculated by using the Euler sequence (see for instance \cite[II.8.20.1]{Ha}). In this case, the Hirzebruch-Milnor class may be more complicated than the Hirzebruch class of $X$.
\bs\bs
\vbox{\centerline{\bf 2. Proofs of Theorem~2 and Proposition~1}
\bsn
In this section we prove Theorem~2 and Proposition~1 by using the short exact sequences associated with nearby and vanishing cycle functors of mixed Hodge modules.}
\msn
{\bf 2.1.~Construction.} In the notation of the introduction, set $S:=\C^{r+1}$ with coordinates $t_1,\dots,t_{r+1}$. We may assume $m\ges 2$ so that $r=\dim\Si_X$, since the case $m=1$ is treated in \cite{MSS1}.
\sk
Let $\Xc\subset Y\times S$ be the hypersurface such that, for any $\ab=(a_1,\dots,a_{r+1})\in S$, we have
$$\Xc\cap(Y\times\{\ab\})=X_{\ab,r+1}\,\bl(:=\{s_{\ab,r+1}=0\}\subset Y{\br)},
\leqno(2.1.1)$$
where $s_{\ab,r+1}$ is as in the introduction (with $n$ replaced by $r+1$ if necessary).
\sk
For $j\in[0,r+1]$, set
$$S_j:=\{t_k=0\,\,(k>j)\}\subset S,\q\Xc_j:=\Xc\times_SS_j\subset\Xc.
\leqno(2.1.2)$$
Note that $\dim S_j=j$, and the fiber of $\Xc_j\to S_j$ over $\ab=(a_1,\dots,a_j)\in S_j$ coincides with $X_{\ab,j}$ in the introduction ($j\in[0,r+1]$). For $j\in[1,r+1]$, set
$$S_{\ab,j}:=\{t_k=a_k\,\,(k<j),\,\,t_k=0\,\,(k>j)\}\subset S_j,\q Y_{\ab,j}:=\Xc\times_SS_{\ab,j}\subset\Xc_j.
\leqno(2.1.3)$$
Then $Y_{\ab,j}$ is a one-parameter family containing $X_{\ab,j-1}$ over $t_j=0$ and $X_{\ab,j}$ over $t_j=a_j$. By the definitions of $U_j$, $f_{\ab,j-1}$ in the introduction, there are natural isomorphisms
$$U_j\cong Y_{\ab,j}\setminus(X'_j\times S_{\ab,j})\,(\subset Y\times S_{\ab,j})\q\q(j\in[1,r+1]),
\leqno(2.1.4)$$
such that $f_{\ab,j-1}$ on the left-hand side is identified with $t_j$ on the right-hand side, since
$$f_{\ab,j-1}=(s_{\ab,j-1}/s_j^{\prime\,m})|_{U_j}=\bl((s-\msum_{k=1}^{j-1}\,a_ks^{\prime\,m}_k)/s^{\prime\,m}_j\br)|_{U_j}=a_j\q\h{on}\q X_{\ab,j}\setminus X'_j.$$
Here the complex number $a_j$ is identified with the variable $t_j$ so that the disjoint union of $X_{\ab,j}$ with $a_j$ varying (including $a_j=0$) is identified with $Y_{\ab,j}$ (where $X_{\ab,j}$ for $a_j=0$ is identified with $X_{\ab,j-1}$.)
\sk
By (2.1.4) we get the isomorphism
$$\varphi_{t_j}\Q_{h,Y_{\ab,j}}=(i_{\Si^{\circ}_{j-1},\Si_{j-1}})_!\varphi_{f_{\ab,j-1}}\Q_{h,U_j}\q\q(j\in[1,r+1]).
\leqno(2.1.5)$$
Indeed, the following is shown in \cite[Section 2.4]{MSS2}:
$$\varphi_{t_j}\Q_{h,Y_{\ab,j}}|_{\Si_{j-1}\cap X'_j}=0.
\leqno(2.1.6)$$
\sk
For $j\in[1,r+1]$, we have the short exact sequences of mixed Hodge modules on $X_{\ab,j}$
$$0\to\Q_{h,X_{\ab,j-1}}[n-1]\to\psi_{t_j}\Q_{h,Y_{\ab,j}}[n-1]\to\varphi_{t_j}\Q_{h,Y_{\ab,j}}[n-1]\to0,
\leqno(2.1.7)$$
since the $a_j$ are sufficiently general.
\msn
{\bf 2.2.~Proof of Theorem~2.} We may assume $m\ges 2$ as in (2.1). There are non-empty Zariski-open subsets $\Sci_j\subset S_j$ ($j\in[1,r+1]$) such that the $D_j:=S_j\setminus \Sci_j$ are divisors on $S_j$ containing $S_{j-1}$ and satisfying
$$D_{j-1}\supset\overline{(D_j\setminus S_{j-1})}\cap S_{j-1}\q(j\in[2,r+1]),
\leqno(2.2.1)$$
and moreover, by setting
$$\Xcc_j:=\Xc\times_S\Sci_j\subset\Xc_j,\q\Ycc_j:=\Xc\times_S\bl(S_j\setminus\overline{(D_j\setminus S_{j-1})}{\br)}\subset\Xc_j,$$
there are short exact sequences of mixed Hodge modules on $\Xcc_{j-1}$ for $j\in[1,r+1]$
$$0\to\Q_{h,\Xcc_{j-1}}[d_{j-1}]\to\psi_{t_j}\Q_{h,\Ycc_j}[d_{j-1}]|_{\Xcc_{j-1}}\to\varphi_{t_j}\Q_{h,\Ycc_j}[d_{j-1}]|_{\Xcc_{j-1}}\to0,
\leqno(2.2.2)$$
where $d_{j-1}:=\dim\Xc_{j-1}\,(=n+j-2)$. We may furthermore assume
$$\Supp\bl(\varphi_{t_j}\Q_{h,\Ycc_j}|_{\Xcc_{j-1}}{\br)}=\Si_{j-1}\times\Sci_{j-1}\subset\Xcc_{j-1},
\leqno(2.2.3)$$
$$\h{$\varphi_{t_j}\Q_{h,\Ycc_j}|_{\Si_{j-1}\times\Sci_{j-1}}\,$ is locally constant over $\Sci_{j-1}$,}
\leqno(2.2.4)$$
by shrinking $\Sci_{j-1}$ if necessary. (Here $\Si_{j-1}=\mcap_{k=1}^{j-1}\,X'_k\cap\Si_X=\Sing\,X_{\ab,j-1}$ since $m\ges 2$.) For the proof of (2.2.4) we use the deformation to the normal bundle in (2.4) below together with a Thom-Sebastiani theorem in Theorem~(3.2) as well as Remarks~(2.5) below.
\sk
These imply by decreasing induction on $k\in[2,j-1]:$
$$\h{$\psi_{t_k}\cdots\psi_{t_{j-1}}(\varphi_{t_j}\Q_{h,\Ycc_j}|_{\Si_{j-1}\times\Sci_{j-1}})|_{\Si_{j-1}\times\Sci_{k-1}}\,$ is locally constant over $\Sci_{k-1}$.}
\leqno(2.2.5)$$
\sk
In the notation of (2.1) we have the isomorphisms
$$\psi_{t_j}\Q_{h,Y_{\ab,j}}=\psi_{t_j}\Q_{h,\Ycc_j}|_{X_{\ab,j-1}},\q\varphi_{t_j}\Q_{h,Y_{\ab,j}}=\varphi_{t_j}\Q_{h,\Ycc_j}|_{X_{\ab,j-1}},
\leqno(2.2.6)$$
such that the restriction of the short exact sequence (2.2.2) to $X_{\ab,j-1}$ is identified with (2.1.7), since the $a_j$ are sufficiently general.
(Indeed, the $V$-filtration induces the $V$-filtration on the restriction to the transversal slice $Y_{\ab,j}\subset\Ycc_j$ passing through $X_{\ab,j-1}\subset\Xcc_{j-1}$, see \cite[Theorem 1.1]{DMST}. Moreover this restriction morphism induces a bistrict surjection for $(F,V)$, see \cite[Lemma 4.2]{DMST}.
So the assertion follows, since the weight filtration $W$ is given by the relative monodromy filtration.)
\sk
We then get the assertion (4) in Theorem~2 by applying the iterations of nearby cycle functors $\psi_{t_k}$ ($k<j$) to (2.2.2) and using (2.1.5), (2.2.6).
Here we also need \cite{Sch4} (or \cite[Proposition 3.3]{MSS1}) together with \cite[Theorem 7.1]{Ve} in order to show that the virtual Hirzebruch class $T^{\,\vir}_{y*}(X)$ can be obtained by applying the iteration of the nearby cycle functors $\psi_{t_j}$ ($j\les r+1$) to $\Q_{h,\Ycc_{r+1}}$.
This finishes the proof of Theorem~2.
\msn
{\bf 2.3.~Proof of Proposition~1.} The limit in the formula (7) is defined by
$$\lim_{a_{j+1}\to 0}T_{y*}(X_{\ab,j+1}):=T_{y*}(\psi_{t_{j+1}}\Q_{h,Y_{\ab,j+1}}),
\leqno(2.3.1)$$
in the notation of (2.1), where the index $j$ is shifted by 1. The assertion (7) then follows from (2.1.5), (2.1.7) and (2.2.4--6).
This finishes the proof of Proposition~1.
\ms
We review here some basics of deformations to normal bundles which will be needed in the proof of Theorem~3.
\msn
{\bf 2.4~Deformations to normal bundles.}
In the notation of the introduction, set
$$x'_i:=(s'_i/s'_j)|_{U_j}\q(i=1,\dots,j-1).$$
Since $Z^{\circ}_{j-1}=\mcap_{i=1}^{j-1}\{x'_i=0\}\subset U_j=Y\setminus X'_j$, we have the decomposition
$$V_j:=N_{Z^{\circ}_{j-1}/U_j}=Z^{\circ}_{j-1}\times\C^{j-1},
\leqno(2.4.1)$$
where the left-hand side is the normal bundle of $Z^{\circ}_{j-1}$ in $U_j$.
The total deformation space $\Uc_j$ of $U_j$ to the normal bundle $V_j$ can be defined by
$$\Uc_j:={\rm Spec}_{U_j}\bl(\mopl_{i\in\Z}\,\I^{\,i}_{Z^{\circ}_{j-1}}\otimes t^{-i}\br),$$
where $\I_{Z^{\circ}_{j-1}}\subset\OO_{U_j}$ is the ideal of $Z^{\circ}_{j-1}$, and $\I^{\,i}_{Z^{\circ}_{j-1}}=\OO_{U_j}$ for $i\les 0$. We can identify $\Uc_j$ with a relative affine open subset of the blow-up of $U_j\times\C$ along $Z^{\circ}_{j-1}\times\{0\}$, on which the functions
$$x_i:=x'_i/t\q(i=1,\dots,j-1)$$
are well-defined, where $t$ is the coordinate of the second factor of $U_j\times\C$, that is, the parameter of deformation.
Note that the normal bundle $V_j$ is contained in $\Uc_j$ as the fiber over $t=0$ (which is a relative affine open subset of the exceptional divisor of the blow-up), since
$$V_j={\rm Spec}_{U_j}\bl(\mopl_{i\ges 0}\,\I^{\,i}_{Z^{\circ}_{j-1}}/\I^{\,i+1}_{Z^{\circ}_{j-1}}\otimes t^{-i}\br).$$
Set
$$z_i:=x_i|_{V_j}\q\q(i=1,\dots,j-1).$$
These give the decomposition (2.4.1) by inducing coordinates of the second factor of (2.4.1).
\sk
Consider now the following function on $\Uc_j:$
$$\ft_{\ab,j-1}:=\pi_j^*(s/s^{\prime\,m}_j)|_{U_j}-\msum_{i=1}^{j-1}\,a_i\,x_i^m,
\leqno(2.4.2)$$
where $\pi_j:\Uc_j\to U_j$ is the canonical morphism.
Restricting over $t=1$, we have
$$\ft_{\ab,j-1}|_{t=1}=f_{\ab,j-1},
\leqno(2.4.3)$$
On the other hand, restricting over $t=0$, we get
$$\ft_{\ab,j-1}|_{V_j}=pr_1^*\bl(f_{\ab,j-1}|_{Z^{\circ}_{j-1}}\br){}-\msum_{i=1}^{j-1}\,a_i\,z_i^m=pr_1^*f'_{j-1}-\msum_{i=1}^{j-1}\,a_i\,z_i^m,
\leqno(2.4.4)$$
where $pr_1:V_j\to Z^{\circ}_{j-1}$ is the canonical projection, and $f'_{j-1}$ is defined just before Theorem~1. To (2.4.4) we can apply a Thom-Sebastiani theorem (see Theorem~(3.2) below).
\msn
{\bf Remarks~2.5.} (i) The restriction of $\ft_{\ab,j-1}$ to a sufficiently small neighborhood $\Uc'_j$ of $Z^{\circ}_{j-1}\times\C$ in $\Uc_j$ is a topologically locally trivial family parametrized by $t\in\C$, since the $s'_j$ are sufficiently general. The vanishing cycle functor along $\ft_{\ab,j-1}$ for $\Q_{h,\Uc'_j}$ commutes with the restriction to $t=c$ for any $c\in\C$. This follows from \cite{DMST} as in the proof of (2.2.6). (Here analytic mixed Hodge modules are used.)
\sk
(ii) Let $\M$ be a mixed Hodge module on $X\times\C$. Assume its underlying $\Q$-complex is isomorphic to the pull-back of a $\Q$-complex on $X$ by the first projection $p:X\times\C\to X$. Then $\M$ is isomorphic to the pull-back of a mixed Hodge module on $X$ by $p$ up to a shift of a complex. Indeed, we have the isomorphism $p^*p_*\M\to\M$ in this case.
(We apply this to $X=Z^{\circ}_{j-1}$ and $\M=\varphi_{\ft_{\ab,j-1}}\Q_{h,\Uc_j}[d_{\Uc_j}-1]|_{Z^{\circ}_{j-1}\times\C}$.)
\msn
{\bf 2.6.~Proof of Corollary~1.} We specialize Theorem~1 to $\y = 1$ (using \cite[Proposition 5.21]{Sch3}) and restrict to the degree 0 part. Since the vanishing cycle functor commutes with the restrictions to the intersections of $X'_j$ (see \cite[Lemma 4.3.4]{Sch2}, \cite{DMST}), it implies that there is $M(X)\in\HH_{\ssb}(\Si_X)$ satisfying $c^{\vir}(X)-c(X)=(i_{\Si_X,X})_*M(X)$, and
$$M(X)=\msum_{j=0}^r\,(1-m)^j\bl((i_{\Si_j,\Si_X})_*c(\EM_X|_{\Si_j})-(i_{\Si_{j+1},\Si_X})_*c(\EM_X|_{\Si_{j+1}})\br).
\leqno(2.6.1)$$
So we get the first equality of (10) in Corollary~1. For the second equality of (10), we have to show
$$(i_{\Si_j,\Si_X})_*c(\EM_X|_{\Si_j})=c(\EM_X)\cap\bl(\tfrac{c_1(L)}{1+c_1(L)}\br)^j.
\leqno(2.6.2)$$
This follows from an analogue of Verdier's Riemann-Roch theorem for the MacPherson Chern classes under transversal pullbacks, see \cite[Corollary 0.1]{Sch1} or \cite[Corollary 2.7]{Sch6}. (Note that the denominator $c_1(L)^j$ in (2.6.2) corresponds to $(i_{\Si_j,\Si_X})_*(i_{\Si_j,\Si_X})^!$.) This finishes the proof of Corollary~1.
\msn
{\bf Remark~2.7.} We can show directly the equality between the last terms of (10) and (12). (This would give a certain reliability to the calculations in (2.1--4).) It is enough to show
$$\tfrac{mx}{1+mx}=\msum_{j\ges 1}\,m(1-m)^{j-1}\bl(\tfrac{x}{1+x}\br)^j\q\h{in}\,\,\,\Z[[x]],
\leqno(2.7.1)$$
since $\tfrac{1}{1+mx}=1-\tfrac{mx}{1+mx}$, and the action of $c_1(L)^{r+1}$ on $\HH_{\ssb}(\Si_X)$ vanishes by $r=\dim\Si_X$. (Note that the cap product with $m\,c_1(L)\,\bl(=c_1(L^{\otimes m})\br)$ corresponds to $(i_{\Si'_X,\Si_X})_*(i_{\Si'_X,\Si_X})^!$.)
\sk
Set $z:=\tfrac{x}{1+x}$ so that $x=\tfrac{z}{1-z}$ in $\Z[[x]]$. Then (2.7.1) holds since
$$\tfrac{mx}{1+mx}=\tfrac{mz}{1-(1-m)z}=\msum_{j\ges 1}\,m(1-m)^{j-1}z^j.
\leqno(2.7.2)$$
\msn
{\bf Example~2.8.} Let $Y_i$ ($i=1,2$) be a smooth hypersurface of degree $a_i\ges 2$ in $\PP^n$ ($n\ges 3$) with $g_i$ a defining polynomial. Assume $Y_1$, $Y_2$ intersect transversally so that $Y_1\cap Y_2$ is smooth. Take $m\in\N$ with $b_i:=m/a_i\in\N$ ($i=1,2$). Let $X$ be a hypersurface of degree $m$ defined by
$$g:=g_1^{b_1}-c\,g_2^{b_2}\q\h{for}\,\,\,\,c\in\C^*\,\,\,\h{sufficiently general,}$$
so that $\Sing\,X=Y_1\cap Y_2$, where we assume $b_i\ges 2$ ($i=1,2$). We have in this case
$$\EM_X(x)=-\mprod_{i=1}^2(b_i-1)\q\h{for any}\,\,\,x\in\Sing\,X=Y_1\cap Y_2.
\leqno(2.8.1)$$
Let $X'$ be a sufficiently general smooth hypersurface of degree $m$ in $\PP^n$. We denote by $\chi(n';a_1,a_2)$ the Euler number of a smooth complete intersection of degree $(a_1,a_2)$ in $\PP^{n'}$ ($n'\les n$), see Remark~(2.9) below. Corollary~2 then implies that $\chi(X')-\chi(X)$ is given by
$$-\mprod_{i=1}^2\bl(\tfrac{m}{\,a_i}-1\br)\cdot\bl(\chi(n;a_1,a_2)-\msum_{j=1}^{n-2}\,m(1-m)^{j-1}\chi(n{-}j;a_1,a_2)\br).
\leqno(2.8.2)$$
On the other hand, (13) implies that $\chi(X')-\chi(X)$ is given by
$$-\mprod_{i=1}^2\bl(\tfrac{m}{\,a_i}-1\br)\cdot\bl(\chi(n;a_1,a_2)-\chi(n;a_1,a_2,m)\br),
\leqno(2.8.3)$$
where $\chi(n;a_1,a_2,m)=\chi(Y_1\cap Y_2\cap X')$. The Euler numbers of smooth complete intersections can be calculated as in Remark~(2.9) below. We have the equality between (2.8.2) and (2.8.3) as a special case of Remark~(2.7). Note that, fixing $a_1,a_2$ (for instance, setting $a_1=a_2=2$), it is enough to calculate only $\chi(n';a_1,a_2)$ in (2.8.2) for {\it any} $m\in\N$ with $m/a_i\in\N$, although we have to calculate also $\chi(n;a_1,a_2,m)$ for every such $m$ in (2.8.3).
\sk
The above argument can be generalized to the case where $g=\prod_{i\les e}\,g_i^{b_i}-c\,\prod_{i>e}\,g_i^{\,b_i}$ with $Y_i:=g_i^{-1}(0)\subset\PP^n$ smooth of degree $a_i$, $\mcup_i\,Y_i$ has normal crossings, and $\sum_{i\les e}a_ib_i=\sum_{i>e}a_ib_i$. Corollary~2 is better when we fix the $g_i$ with $b_i$ varying.
\msn
{\bf Remark~2.9.} Let $X$ be a smooth complete intersection of degree $(d_1,\dots,d_k)$ in $\PP^n$, that is, $X$ is defined by homogeneous polynomials $g_i$ of degree $d_i$ ($i\in[1,k]$) and $\dim X=n-k$. As a special case of Hirzebruch theory \cite{Hi}, it is well-known (see for instance \cite{Nav}) that the Euler number of $X$, denoted by $\chi(n;d_1,\dots,d_k)$, is given by the coefficient of $x^n$ in
$$\gamma(n;d_1,\dots,d_k):=(1+x)^{n+1}\,\mprod_{i=1}^k\,\tfrac{d_ix}{1+d_ix}\in\Z[x]/(x^{n+1})\,\bl(=H^{\ssb}(\PP^n)\br).
\leqno(2.9.1)$$
More precisely, $\gamma(n;d_1,\dots,d_k)$ expresses the {\it image} of the (virtual) Chern class of $X$ under the inclusion $X\into\PP^n$ using Poincar\'e duality (see \cite[1.4.5]{MSS2}), and the product of the $d_ix$ in the numerator comes from this. (Note that the intersection number of $X$ with a linear subspace in $\PP^n$ generating $H^{2n-2k}(\PP^n,\Z)$ is $d_1\cdots d_k$.) For instance, if $n=2$, $k=1$, then the Euler number of a plane curve of degree $d$ is given by $3d-d^{\,2}=2-(d-1)(d-2)$ as is well-known.
\bs\bs
\vbox{\centerline{\bf 3. Application of Thom-Sebastiani theorem}
\bsn
In this section we prove Theorems~3 and 4 after explaining the Thom-Sebastiani theorem for underlying filtered $\D$-modules of constant mixed Hodge modules.}
\msn
{\bf 3.1.~Algebraic microlocalization.} Let $Y$ be a smooth complex algebraic variety (or a connected complex manifold) with $f$ a non-constant function on $Y$, that is $f\in\Gamma(Y,\OO_Y)\setminus\C$. Set
$$(\B_f,F):=(i_f)_*^{\D}(\OO_Y,F)=(\OO_Y[\dd_{t}],F[-d_Y]),
\leqno(3.1.1)$$
with $d_Y:=\dim Y$, and $F[m]_p=F[m]^{-p}=F^{m-p}=F_{p-m}$ in general. Here
$(i_f)_*^{\D}$ is the direct image of filtered $\D$-modules with $i_f:Y\into Y\times\C$ the graph embedding by $f$, and $t$ is the coordinate of the second factor of $Y\times\C$, see (1.2.4).
The last isomorphism is as filtered $\OO_Y[\dd_{t}]$-modules, where the sheaf-theoretic direct image $(i_f)_*$ is omitted to simplify the notation. Note that $g\,\dd_t^j\in\OO_Y[\dd_{t}]$ is actually identified with $g\,\dd_t^j\delta(t-f)$, where
$$\delta(t-f):=\h{$\frac{1}{t-f}$}\in\OO_{Y\times\C}\bl[\h{$\frac{1}{t-f}$}\br]/\OO_{Y\times\C}.$$
Moreover $\delta(t-f)$ is also identified with $f^s$ (and $-\dd_tt$ with $s$), see for instance \cite{Mal}.
\sk
In this section, the Hodge filtration $F$ is indexed as in the case of {\it right} $\D$-modules (since there is a shift of filtration under the direct images by closed embeddings for left $\D$-modules, see (1.2.4)). So we have
$$\Gr_p^F\B_f=\begin{cases}\OO_Y\dd_{t}^{\,p+d_Y}&\h{if}\,\,\,\ p\ges-d_Y,\\
\,0&\h{otherwise}.\end{cases}
\leqno(3.1.2)$$
By (1.2.2) this does not cause a problem when we use the de Rham functor $\DR$.
\sk
Let $\BB_f$ be the {\it algebraic microlocalization} of $\B_f$ (see \cite{mic}), that is,
$$(\BB_f,F)=(\OO_Y[\dd_{t},\dd_{t}^{-1}],F)\q\h{with}\q\Gr_p^F\BB_f=\OO_Y\dd_{t}^{\,p+d_Y}\q(p\in\Z).
\leqno(3.1.3)$$
Let $V$ be the microlocal $V$-filtration on $\BB_f$ along $t=0$, see \cite{mic}. This is defined by modifying the $V$-filtration of Kashiwara \cite{Ka} and Malgrange \cite{Mal} on $\B_f$, see also \cite{MSS3}.
\sk
By construction there are canonical isomorphisms
$$\can:\Gr_V^{\al}(\B_f,F)\simto\Gr_V^{\al}(\BB_f,F)\q\q(\al<1),
\leqno(3.1.4)$$
$$\dd_t^{\,k}:(\BB_f;F,V)\simto(\BB_f;F[-k],V[-k])\q\q(k\in\Z),
\leqno(3.1.5)$$
where $\,\can\,$ in (3.1.4) for $\al=1$ is strictly surjective by \cite[Lemma 5.1.4 and Proposition 5.1.14]{mhp}. Indeed, setting
$$X:=f^{-1}(0)\subset Y,$$
we have a short exact sequence of mixed Hodge modules
$$0\to\Q_{h,X}[d_X]\to\psi_{f,1}\Q_{h,Y}[d_X]\buildrel{\can\,\,}\over\longrightarrow\varphi_{f,1}\Q_{h,Y}[d_X]\to 0,
\leqno(3.1.6)$$
and $\,\can\,$ in (3.1.4) for $\al=1$ is identified with the underlying morphism of filtered $\D_Y$-modules of the morphism $\,\can\,$ in (3.1.6).
\sk
Note that we have for $\al\in(-1,0]$ (see for instance \cite[1.1.9]{MSS3})
$$\Gr^F_p\Gr_V^{\al}\BB_f=0\,\,\,\,(p<-d_X).
\leqno(3.1.7)$$
\sk
We have a Thom-Sebastiani theorem as below. This is a special case of the assertion mentioned in \cite[Remark 4.5]{mic}, and follows from \cite[Theorem 1.2]{MSS3} or \cite{ts} (see also \cite{DL}, \cite{GLM} for the motivic version, and \cite{ScSt}, \cite{Va} for the isolated hypersurface singularity case).
\msn
{\bf Theorem~3.2.} {\it Let $Y_a$ be a smooth complex algebraic variety $($or a connected complex manifold$)$ with $f_a$ a non-constant function, that is, $f_a\in\Gamma(Y_a,\OO_{Y_a})\setminus\C$, for $a=1,2$. Set $Y=Y_1\times Y_2$ with $f=f_1+f_2$. Then there are canonical isomorphisms of filtered $\D_Y$-modules for $\al\in(-1,0]:$
$$\aligned\Gr^{\al}_V(\BB_f,F)&=\mopl_{\al_1\in I(\al)}\,\Gr^{\al_1}_V(\BB_{f_1},F)\boxtimes\Gr^{\al-\al_1}_V(\BB_{f_2},F)\\
&\q\oplus\,\mopl_{\al_1\in J(\al)}\,\Gr^{\al_1}_V(\BB_{f_1},F)\boxtimes\Gr^{\al-1-\al_1}_V(\BB_{f_2},F[-1]),\endaligned
\leqno(3.2.1)$$
by replacing $Y_a$ with an open neighborhood of $X_a:=f_a^{-1}(0)$ in $Y_a$ {$(a=1,2)$} if necessary, where}
$$I(\al):=(-1,0]\cap[\al,\al+1)\q\q J(\al):=(-1,0]\cap[\al-1,\al).$$
\msn
{\bf Note.} Setting $\al_2=\al-\al_1$, $\,\al'_2=\al-1-\al_1$, we have
$$\al_1\in I(\al)\iff\al_1,\,\al_2\in(-1,0],\q\q\al_1\in J(\al)\iff\al_1,\,\al'_2\in(-1,0].
\leqno(3.2.2)$$
\msn
{\bf Remarks~3.3.} (i) Since the de Rham functor $\DR$ is compatible with the exterior product $\boxtimes$, we can deduce from Theorem~(3.2) the following isomorphisms of complexes of $\OO_Y$-modules for $\al\in(-1,0]$, $p\in\Z$:
$$\aligned&\Gr_F^p\DR_Y(\Gr_V^{\al}\BB_f)\\
&=\mopl_{\al_1+\al_2=\al,\,p_1+p_2=p}\,\Gr_F^{p_1}\DR_{Y_1}(\Gr_V^{\al_1}\BB_{f_1})\boxtimes\Gr_F^{p_2}\DR_{Y_2}(\Gr_V^{\al_2}\BB_{f_2})\\
&\q\oplus\,\mopl_{\al_1+\al_2=\al-1,\,p_1+p_2+1=p}\,\Gr_F^{p_1}\DR_{Y_1}(\Gr_V^{\al_1}\BB_{f_1})\boxtimes\Gr_F^{p_2}\DR_{Y_2}(\Gr_V^{\al_2}\BB_{f_2}),\endaligned
\leqno(3.3.1)$$
where $\al_1,\al_2\in(-1,0]$, and $F^p=F_{-p}$. (Note that $\Gr_F^{p_2}=\Gr_{F[-1]}^{p'_2}$ with $p'_2:=p_2+1$.)
\sk
Setting $\la_a=\exp(-2\pi i\al_a)$, we have $\al_a=-\ell(\la_a)$ ($a=1,2$), and
$$\al_1+\al_2\les-1\iff\ell(\la_1)+\ell(\la_2)\ges 1.
\leqno(3.3.2)$$
If these equivalent conditions are satisfied, then (3.3.1) says that the index $p$ of the Hodge filtration $F$ increases by $1$.
This is very important for the proofs of Theorems~3 and 4.
\sk
(ii) In the notation of (3.1), assume $Y_2=\C$ with coordinate $x_2$, and $f_2=a\,x_2^m$ with $a\in\C^*$. Then the Milnor fiber $F_{\!f_2,0}$ consists of $m$ points, and
$$\Ht^0(F_{\!f_2,0},\C)_{\la}=\begin{cases}\C&\h{if}\,\,\,\la^m=1\,\,\,\h{and}\,\,\,\la\ne 1,\\
\,0&\h{otherwise,}\end{cases}
\leqno(3.3.3)$$
where $\Ht$ denotes the reduced cohomology.
This implies
$$\Gr_F^{p_2}\DR_{Y_2}(\Gr_V^{\al_2}\BB_{f_2})\cong\begin{cases}\C&\h{if}\,\,\,p_2=0,\,\,\al_2\in\bl\{\frac{1}{m},\dots,\frac{m-1}{m}\br\},\\
\,0&\h{otherwise.}\end{cases}
\leqno(3.3.4)$$
Note that we have in this case
$$\varphi_f\C_Y=\varphi_{f_1}\C_{Y_a}\otimes_{\C}\Ht^0(F_{\!f_2,0},\C)[-1].
\leqno(3.3.5)$$
\msn
{\bf 3.4.~Proofs of Theorems~3 and 4.} We get (5) in Theorem~3 by Proposition~(1.5), and the assertion (6) follows from (2.4.3--4), (3.3.1--2), (3.3.4--6) together with Remarks~(2.5).
Theorem~4 also follows by using (3.3.1--2) and (3.3.6) together with the compatibility of the Todd class transformation $td_*:K_0(X)\to\HH_{\ssb}(X)$ with cross products (or K\"unneth maps), see \cite[Section III.3]{BFM}.
This finishes the proofs of Theorems~3 and 4.
\msn
{\bf Remarks~3.5.} (i) By the proof of Theorem~4, the assertion holds at the level of Grothendieck groups. We have the following equality in $K_0(\Si_X)\bl[\y^{1/e}\br]:$
$$\DR_{\y}[\varphi_f\Q_{h,Y},T_s]=-\DR_{\y}[\varphi_{f_1}\Q_{h,Y_1},T_s]\,\boxtimes\,\DR_{\y}[\varphi_{f_2}\Q_{h,Y_2},T_s].
\leqno(3.5.1)$$
\sk
(ii) Theorem~4 does not necessarily hold if there are non-zero critical values $c_a$ of $f_a$ ($a=1,2$) with $c_1+c_2=0$, since the last condition is equivalent to the condition $\Si_X\ne\Si_{X_1}\times\Si_{X_2}$.
However, we can apply Theorem~4 with $f_a$ replaced by $f_a-c_a$ ($a=1,2$) in the above case.
\sk
(iii) In Theorem~4, $X$ is never compact even if $X_1,X_2$ are compact, since $X$ is the inverse image of the anti-diagonal of $\C\times\C$ by $f_1\times f_2$.
Extending the situation in Theorem~4, we may consider the case where $f_1$, $f_2$ are proper morphisms from smooth varieties to $\PP^1$ and $X$ is defined by the inverse image of the anti-diagonal of $\PP^1\times\PP^1$. In this case $X$ is compact. We can apply Theorem~4 by choosing an appropriate local coordinate of $\PP^1$ on a neighborhood of each critical value $c_a$ of $f_a$ ($a=1,2$) with $(c_1,c_2)$ belonging to the anti-diagonal of $\PP^1\times\PP^1$.
\bs\bs
\vbox{\centerline{\bf 4. Relation with rational and Du Bois singularities}
\bsn
In this section some relations with rational and Du Bois singularities are explained.}
\msn
{\bf 4.1.~Primitive decomposition of nearby cycles.} Let $Y$ be a smooth complex algebraic variety (or a connected complex manifold), and $f$ be a non-constant function on $Y$, that is, $f\in\Gamma(Y,\OO_Y)\setminus\C$.
Assume $X:=f^{-1}(0)\subset Y$ is {\it reduced\,} in this section.
\sk
With the notation of (3.1), we first show the short exact sequence
$$0\to\omt_X{\otimes_{\OO_X}}\om_X^{\vee}\to F_{-d_X}\Gr_V^1\B_f\buildrel\can\over\longrightarrow F_{-d_X}\Gr_V^1\BB_f\to 0.
\leqno(4.1.1)$$
Here
$$\omt_X:=(\rho)_*\om_{\X}\subset\om_X,\q\om_X^{\vee}:=\Hc om_{\OO_X}(\om_X,\OO_X),$$
with $\rho:\X\to X$ a resolution of singularities, and
$$\om_X=\om_Y\otimes_{\OO_Y}\OO_X,$$
since $X$ is globally defined by $f$, see \cite[Lemma 2.9]{rat}.
\sk
For the proof of (4.1.1), we use the monodromy filtration $W$ on $\Gr_V^1\B_f$ shifted by $d_X$, which is uniquely characterized by the following conditions:
$$\aligned N(W_i\,\Gr_V^1\B_f)&\subset W_{i-2}\,\Gr_V^1\B_f\q(i\in\Z),\\ N^i:\Gr_{d_X+i}^W\,\Gr_V^1\B_f&\simto\Gr^W_{d_X-i}\,\Gr_V^1\B_f\q(i>0).\endaligned
\leqno(4.1.2)$$
The {\it primitive part} is defined by
$${}^P\Gr^W_{d_X+i}\,\Gr_V^1\B_f:={\rm Ker}\,N^{i+1}\subset\Gr^W_{d_X+i}\,\Gr_V^1\B_f\q(i\ges 0),
\leqno(4.1.3)$$
with ${}^P\Gr^W_{d_X+i}\,\Gr_V^1\B_f=0$ ($i<0$).
This implies the {\it primitive decomposition}
$$\Gr^W_j\,\Gr_V^1\B_f=\mopl_{k\ges 0}\,N^k\bl({}^P\Gr^W_{j+2k}\,\Gr_V^1\B_f{\br)}\q(j\in\Z),
\leqno(4.1.4)$$
and the {\it co-primitive part} can be expressed by
$$\aligned N^i\bl({}^P\Gr^W_{d_X+i}\,\Gr_V^1\B_f\br)&={\rm Ker}\bl(\Gr^W_{d_X-i}\,\Gr_V^1\B_f\buildrel{N}\over\longrightarrow\Gr^W_{d_X-i-2}\,\Gr_V^1\B_f{\br)}\\&={\rm Ker}\bl(\Gr^W_{d_X-i}\,\Gr_V^1\B_f\buildrel{\rm can\,\,}\over\longrightarrow\Gr^W_{d_X-i}\,\Gr_V^1\BB_f{\br)}\q(i\ges0).\endaligned
\leqno(4.1.5)$$
Indeed, the first isomorphism follows from the primitive decomposition (4.1.4). As for the last isomorphism, note that the canonical morphism
$$\can:\Gr_V^1\B_f\to\Gr_V^1\BB_f$$
is identified with the morphism
$$\can:\Gr_V^1\B_f\to\Gr_V^0\B_f,$$
which is defined by $-\Gr_V\dd_t$. Moreover, for the latter, we have
$$N={\rm Var}\ssc\can,$$
where $\,{\rm Var}\,$ is defined by $\Gr_Vt$, and is injective, see \cite[5.1.3.4]{mhp}. So the last isomorphism of (4.1.5) also follows.
\sk
Returning to the proof of (4.1.1), we get by (3.1.7)
$$F_{-d_X}\bl(N^i\bl({}^P\Gr^W_{d_X+i}\,\Gr_V^1\B_f\br){\br)}=F_{-d_X-i}\bl({}^P\Gr^W_{d_X+i}\,\Gr_V^1\B_f{\br)}=0\q(i>0),
\leqno(4.1.6)$$
and it follows from \cite[Proposition 2.7]{rat} that
$$F_{-d_X}\bl({}^P\Gr^W_{d_X}\,\Gr_V^1\B_f{\br)}=\omt_X{\otimes_{\OO_X}}\om_X^{\vee},
\leqno(4.1.7)$$
since $\B_f$ is a left $\D$-module. So (4.1.1) follows.
\sk
We denote by ${\rm IC}_X\Q_h$ the mixed Hodge module of weight $d_X$ such that its underlying $\Q$-complex is the intersection complex ${\rm IC}_X\Q$. We have
$${\rm IC}_X\Q_h=\Gr_{d_X}^W(\Q_{h,X}[d_X])={}^P\Gr^W_{d_X}\psi_{f,1}\Q_{h,Y}[d_X],
\leqno(4.1.8)$$
where the two isomorphisms follow from \cite[4.5.9]{mhm} and (3.1.6) together with the primitive decomposition as in (4.1.4).
\sk
Comparing (4.1.1) with (3.1.6), we see that the exactness of (4.1.1) is essentially equivalent to
$$F_{-d_X}(\Q_{h,X}[d_X])=F_{-d_X}({\rm IC}_X\Q_h)=\omt_X.
\leqno(4.1.9)$$
Here we set in general
$$F_{p_0}\M:=F_{p_0}M,
\leqno(4.1.10)$$
if $(M,F)$ is the underlying filtered {\it right} $\D$-module of a mixed Hodge module $\M$, where
$$p_0:=\min\bl\{p\in\Z\mid\Gr_p^FM\ne 0\br\}.$$
These are independent of embeddings of algebraic varieties into smooth varieties as long as {\it right} $\D$-modules are used.
\msn
{\bf 4.2.~Rational singularities.} With the notation and assumption of (4.1), we have the following canonical isomorphism by (4.1.1):
$$(\om_X/\omt_X){\otimes_{\OO_X}}\om_X^{\vee}=F_{-d_X}(\BB_f/V^{>1}\BB_f),
\leqno(4.2.1)$$
(see also \cite[Theorem 0.6]{rat}), since
$$F_{-d_X}(\B_f/V^{>1}\B_f)=\OO_X,
\leqno(4.2.2)$$
where $V^{>\al}:=V^{\al+\ep}\,\,\,(0<\ep\ll 1$) for $\al\in\Q$ in general.
\sk
Indeed, by (3.1.7), (3.1.4--5), and \cite[3.2.1.2]{mhp}, we have
$$F_{-d_X}V^{>0}\B_f=F_{-d_X}\B_f=\OO_Y,\q F_{-d_X}V^{>1}\B_f=t(F_{-d_X}V^{>0}\B_f).
\leqno(4.2.3)$$
Thus (4.2.2) and (4.2.1) follow.
\sk
As a corollary of (4.2.1), we see that $X$ has only {\it rational singularities} if and only if
$$F_{-d_X}(\BB_f/V^{>1}\BB_f)=0,\q\h{or equivalently}\q F_{-d_X}(\varphi_f\Q_{h,Y}[d_X])=0,
\leqno(4.2.4)$$
under the notation (4.1.10), see also \cite[Theorem 0.6]{rat} (and \cite{exg} in the isolated singularity case). Consider the classes
$$\aligned\bl[(\om_X/\omt_X){\otimes_{\OO_X}}\om_X^{\vee}\br]&=\bl[F_{-d_X}(\BB_f/V^{>1}\BB_f)\br],\\ \bl[\om_X/\omt_X\br]&=\bl[F_{-d_X}(\varphi_f\Q_{h,Y}[d_X])\br]\q\h{in}\,\,\,K_0(\Si_X).\endaligned
\leqno(4.2.5)$$
These may be called the {\it irrationality} of the singularities of $X$ at least in the $\Si_X$ projective case by the argument after (4.2.7) below.
In the isolated singularity case, its dimension is called the geometric genus, see for instance \cite{exg}.
\sk
If $X$ has only rational singularities, then the last condition in (4.2.4) implies
$$\DR_y\bl[\varphi_f\Q_{h,Y}[d_X]\br]\big|{}_{y^{d_X}}=\bl[F_{-d_X}(\varphi_f\Q_{h,Y}[d_X])\br]=0\q\h{in}\,\,\,\,K_0(\Si_X),
\leqno(4.2.6)$$
(see also (1.2.6)), where $|_{y^{d_X}}$ means taking the coefficient of $y^{d_X}$.
(Recall that {\it right} $\D$-modules are used in (4.1.10).) Moreover we have the following.
$$\h{The converse holds if $\Si_X$ is a {\it projective} variety.}
\leqno(4.2.7)$$
Indeed, if the singularities of $X$ are irrational, then we can show the non-vanishing of (4.2.6) in $K_0(\PP^N)_{\Q}$ with $\PP^N$ projective space containing $\Si_X$ by using the topological filtration in (1.6) together with the positivity (see (1.6.4--5)) of the image by the cycle class map (1.6.3) of the coherent sheaf
$$\F:=F_{-d_X}(\varphi_f\Q_{h,Y}[d_X]).$$
\msn
{\bf 4.3.~Du Bois singularities.} With the notation and assumption of (4.1), let $\DD_{h,X}$ be the dual of $\Q_{h,X}$. Since $\Q_{h,X}[d_X]$ is a mixed Hodge module, so is $\DD_{h,X}[-d_X]$. Then we have the short exact sequence of mixed Hodge modules
$$0\to\varphi_{f,1}\Q_{h,Y}(1)[d_X]\buildrel{\rm Var\,\,}\over\longrightarrow\psi_{f,1}\Q_{h,Y}[d_X]\to\DD_{h,X}(-d_X)[-d_X]\to 0,
\leqno(4.3.1)$$
which is the dual of (3.1.6) (up to a sign). Indeed, $\,{\rm Var}\,$ is the dual of $\,\can\,$ in (3.1.6) up to a sign, see \cite[Section 5.2]{mhp}. The underlying exact sequence of filtered $\D$-modules of (4.3.1) is identified with
$$0\to({\rm Im}\,N,F)\to(\Gr_V^1\B_f,F)\to({\rm Coker}\,N,F)\to 0,
\leqno(4.3.2)$$
and the primitive decomposition (4.1.4) implies
$$\Gr^W_{d_X+i}({\rm Coker}\,N,F)={}^P\Gr^W_{d_X+i}(\Gr_V^1\B_f,F)\q(i\ges 0),
\leqno(4.3.3)$$
since the graded quotients $\Gr^W$ commute with taking the cokernel of $N$, see \cite[Proposition~5.1.14]{mhp}. Setting
$$\omt'_X:=F_0\bl(\DD_{h,X}[-d_X]{\br)}=\Gr_F^0\bl(\DR(\DD_{h,X}[-d_X])\br),$$
we then get by (4.3.1--3) and (4.1.6)
$$\omt'_X=F_{-d_X}\bl(\DD_{h,X}(-d_X)[-d_X]{\br)}=F_{-d_X}\bl(\psi_{f,1}\Q_{h,Y}[d_X]\br).
\leqno(4.3.4)$$
Combined with (4.2.2), these imply
$$(\om_X/\omt'_X){\otimes_{\OO_X}}\om_X^{\vee}=F_{-d_X}(\B_f/V^1\B_f)=F_{-d_X}(\BB_f/V^1\BB_f).
\leqno(4.3.5)$$
Since the dual functor $\DD$ commute with $\DR$ (or rather $\DR^{-1}$, see \cite[Section 2.4.11]{mhp}) and also with $\Gr_F^0$ by definition, we have by the definition of $\omt'_X$ just before (4.3.4)
$$\DD(\omt'_X)=\Gr_F^0\DR(\Q_{h,X}[d_X]).
\leqno(4.3.6)$$
Here the left-hand side is the Grothendieck dual of the $\OO_X$-module $\omt'_X$, and we have
$$\DD(\F):=\R\Hc om_{\OO_X}(\F,\om_X[d_X])\q\h{for}\q\F\in D^b_{\rm coh}(\OO_X).$$
\sk
It follows from (4.3.6) that $X$ has only {\it Du Bois singularities} (see \cite{St3}) if and only if
$$\omt'_X=\om_X,
\leqno(4.3.7)$$
(since $\DD(\om_X)=\OO_X[d_X]$).
This condition is equivalent to the vanishing of the $\OO_X$-modules in (4.3.5). Moreover the last condition is equivalent to
$$F_{-d_X}(\psi_{f,\ne1}\Q_{h,Y}[d_X])=F_{-d_X}(\varphi_{f,\ne1}\Q_{h,Y}[d_X])=0.
\leqno(4.3.8)$$
This implies by using \cite[Theorem 0.1]{BS1}
$$\h{$X$ has only Du Bois singularities if and only if $\lct(f)=1$,}
\leqno(4.3.9)$$
where the {\it log canonical threshold} $\lct(f)$ is defined to be the minimal jumping coefficients as in the introduction. Indeed, we have by \cite[Theorem 0.1]{BS1}
$$\G(\al X)=\Gr_V^{\al}\OO_Y\subset\Gr_V^{\al}(\D_Y[s]f^s)\subset\Gr_V^{\al}\B_f\q\q(\forall\,\al\in\Q),
\leqno(4.3.10)$$
where $\G(\al X):=\J((\al-\ep)X/\J(\al X)$ for $0<\ep\ll 1$, and $\D_Y[s]f^s$ is as in Remark~(4.4)(i) below.
\msn
{\bf Remarks~4.4.} (i) The assertion (4.3.9) is equivalent to \cite[Theorem 0.5]{fil} where the statement is given in terms of the maximal root $-\al_f$ of the Bernstein-Sato polynomial $b_f(s)$. Indeed, it is well-known that
$$\lct(f)=\al_f.
\leqno(4.4.1)$$
This follows, for instance, from \cite[Theorem 0.1]{BS1} (see also (4.3.10)) combined with an assertion in \cite{Mal} (more precisely, the roots of the Bernstein-Sato polynomial consist of rational numbers $-\al$ with $\Gr_V^{\al}(M_f/tM_f)\ne 0$ where $M_f:=\D_Y[s]f^s\subset\B_f$, see also a remark after (3.1.1)). This well-known assertion, however, does not seem to be quoted in \cite{KoSc}, although the theorem in \cite{fil} explained above is mentioned after \cite[Corollary 6.6]{KoSc}, where it is shown that a reduced hypersurface $X\subset Y$ has only Du Bois singularities if and only if $(Y,X)$ is a {\it log canonical pair}, see also Remark~(ii) below.
\sk
(ii) It is well-known (and is easy to show) that $(Y,X)$ is a log canonical pair with $X$ reduced if and only if $\lct(f)=1$. Indeed, let $\rho:(\Yt,\X)\to(Y,X)$ be an embedded resolution. We have
$$\X=\rho^*X=\X'+\msum_i\,m_iE_i,\q\om_{\Yt}=(\rho^*\om_Y)\bl(\msum_i\,\nu_iE_i{\br)}\q(m_i,\nu_i\in\Z_{>0}),
\leqno(4.4.2)$$
where the $E_i$ are the exceptional divisors of $\rho$, and $\X'$ is the proper transform of $X$. (The last equality is equivalent to that ${\rm div}\bl({\rm Jac}(\rho){\br)}=\msum_i\,\nu_iE_i$, where ${\rm Jac}(\rho)$ is the Jacobian of $\rho$ with respect to some local coordinates of $\Yt$, $Y$.) By (4.4.2) we then get
$$\om_{\Yt}(\X')=\bl(\rho^*\om_Y(X){\br)}\bl(\msum_i\,(\nu_i-m_i)E_i\br).
\leqno(4.4.3)$$
Since $\lct(f)=\min\JC(f)$ by definition, the equality (4.4.3) implies
$$\aligned&(Y,X)\,\,\,\h{is a log canonical pair}\iff\nu_i-m_i\ges -1\,\,\,(\forall\,i)\\&\iff (\nu_i+1)/m_i\ges 1\,\,\,(\forall\,i)\iff\lct(f)=1,\endaligned
\leqno(4.4.4)$$
where the first equivalence is by the definition of canonical pairs together with (4.4.3), see \cite{KoSc}. The last equivalence follows from the well-known assertion:
$$\lct(f)=\min\{(\nu_i+1)/m_i\}\,\,\,\,\h{if}\,\,\,\,\lct(f)<1\,\,\,\,\h{or}\,\,\,\min\{(\nu_i+1)/m_i\}<1.
\leqno(4.4.5)$$
By the definition of $\JC(f)$ (see for instance \cite[Section 2.1]{MSS3}), the last assertion can be verified by calculating the integration of $\rho^*(|f|^{-2\al}\om\wedge\overline{\om})$ on $\Yt$, where $\om$ is a nowhere vanishing holomorphic form of degree $d_Y$ locally defined on $Y$, and $\al\in(0,1]$. Indeed, this can be reduced to a well-known assertion saying that we have for $\beta\in\R$, $c\in\R_{>0}$
$$\int_0^cr^{\beta}\ddd r<\infty\iff\beta>-1.$$
\msn
{\bf 4.5.~Proof of Theorem~5.} We first show the second case where $X$ is globally defined by a function $f$ on $Y$.
By the induced polarization on the nearby and vanishing cycle mixed Hodge modules (see \cite[Section 5.2]{mhp}, \cite{dual}), we have the self-dualities
$$\aligned\DD(\varphi_{f,\ne 1}\Q_{h,Y}[d_X])&=(\varphi_{f,\ne 1}\Q_{h,Y}[d_X])(d_X),\\ \DD(\varphi_{f,1}\Q_{h,Y}[d_X])&=(\varphi_{f,1}\Q_{h,Y}[d_X])(d_X+1),\endaligned
\leqno(4.5.1)$$
since $\psi_{f,\ne 1}=\varphi_{f,\ne 1}$, and $d_Y=d_X+1$. (Note that the monodromy filtration is self-dual.) These imply in the notation of (4.1.10)
$$\aligned\DD\bl(\Gr_F^0\DR(\varphi_{f,\ne 1}\Q_{h,Y}[d_X])\br)&=\Gr_F^{d_X}\DR(\varphi_{f,\ne 1}\Q_{h,Y}[d_X])\\&=F_{-d_X}(\varphi_{f,\ne1}\Q_{h,Y}[d_X]),\\ \DD\bl(\Gr_F^0\DR(\varphi_{f,1}\Q_{h,Y}[d_X])\br)&=\Gr_F^{d_X+1}\DR(\varphi_{f,1}\Q_{h,Y}[d_X])\\&=0.\endaligned
\leqno(4.5.2)$$
Recall that {\it right} $\D$-modules are used in (4.1.10).
\sk
So the first assertion of Theorem~5 in the second case follows from the assertion concerning (4.3.8), since $\DD^2=id$ and
$$M_0(X)=td_*\bl[\Gr_F^0\DR(\varphi_f\Q_{h,Y})\br].
\leqno(4.5.3)$$
\sk
To show the converse, assume that $X$ is {\it not} Du Bois, that is,
$$\F:=F_{-d_X}(\varphi_{f,\ne1}\Q_{h,Y}[d_X])\ne 0.
\leqno(4.5.4)$$
By (4.5.2--3) we have to show
$$td_*\bl[\DD(\F)\br]\ne 0\q\h{in}\,\,\,\HH_{\ssb}(\PP^N),
\leqno(4.5.5)$$
where we can replace $\HH_{\ssb}(\Si_X)$ with $\HH_{\ssb}(\PP^N)$ by the compatibility of $td_*$ with the pushforward by proper morphisms.
Then the assertion follows by using the topological filtration on $K_0(\PP^N)_{\Q}$ and $\HH_{\ssb}(\PP^N)$ in (1.6) together with the positivities in (1.6.4--5) (see also an argument after (4.2.7)).
This finishes the proof of Theorem~5 in the second case.
\sk
For the proof in the first case, note that the support of $\F$ in (4.5.4) is independent of the choice of a local defining function of $X$ (where the ambiguity comes from the multiplication by a nowhere vanishing function).
However, we have to take here the direct image $(i_{\Si^{\circ}_0,\Si_X})_!$ of a mixed Hodge module. This can be calculated as in the proof of Proposition~(1.4), and the latter shows that it is enough to take the {\it closure} of the support of the coherent sheaf which gives the non-Du Bois locus. This closure is independent of the choice of $s'_1$, and taking the direct image does not cause a problem as long as $s'_1$ is sufficiently general so that $X'_1=s_1^{\prime\,-1}(0)$ does not contain this support. So the assertion follows. Here it is enough to consider the summand in the formula (4) in Theorem~2 with $j=1$ by using the addition theorem for the log canonical threshold (see \cite[Corollary 1]{MSS3}) together with Theorem~3.
This finishes the proof of Theorem~5.
\msn
{\bf Remark~4.6.} We do not know a priori the support of the coherent sheaf in the above argument, and there might be some problem about the genericity condition on $s'_1$ (that is, the condition that $s_1^{\prime\,-1}(0)$ does not contain the support). It may be better to argue as follows:
\sk
On a dense Zariski-open subset $U$ of the parameter space of $s'_1$, $X'_1=s_1^{\prime\,-1}(0)$ intersects the strata of a Whitney stratification of $X$ transversally so that $M_0(X)$ can be defined. Moreover $M_0(X)$ in the graded pieces of the topological filtration in (1.6) is independent of the choice of $s'$, since it is given by the cycle map, see (1.6). (Here \cite{DMST} is also used.) There is another dense Zariski-open subset $U'$ of the parameter space of $s'_1$ such that $X'_1=s_1^{\prime\,-1}(0)$ does not contain the non-Du Bois locus. We have $U\subset U'$, since the image of the cycle map would vanish if $s'_1\notin U'$. So no problem occurs.
\msn
{\bf 4.7.~Proof of Proposition~2.} The duality isomorphisms in (4.5.1) are compatible with the action of the semisimple part of the monodromy $T_s$, where the $\e(\al)$-eigenspace is the dual of the $\e(-\al)$-eigenspace, see also \cite[2.4.3]{mic}. The argument is then essentially the same as in the proof of Theorem~5 by using the topological filtration in (1.6) together with \cite[Theorem 0.1]{BS1} (see also (4.3.10)) which gives the relation with the jumping coefficients. This finishes the proof of Proposition~2.
\msn
{\bf 4.8.~Isolated hypersurface singularity case.} Let $f:(Y,0)\to(\Delta,0)$ be a germ of a holomorphic function on a complex manifold $Y$ such that $X:=f^{-1}(0)$ has an isolated singularity at $0$, where $\Delta\subset\C$ is an open disk. Let $\mu_f$ be the Milnor number of $f$. As in \cite{St2}, the {\it spectrum}
$${\rm Sp}(f)=\msum_{i=1}^{\mu_f}\,t^{\,\al_{f,i}}\in\Z\bl[t^{1/e}\br]$$
with $\al_{f,i}\les\al_{f,i+1}\,\,\,(i\in[1,\mu_f-1])$ is defined by
$$\#\bl\{i\mid\al_{f,i}=\al\br\}=\dim\Gr_F^p\Ht^{d_Y-1}(F_{\!f,0},\C)_{\e(-\al)}\q\q(p:=[d_Y-\al],\,\,\al\in\Q),
\leqno(4.8.1)$$
where $F_{\!f,0}$ denotes the Milnor fiber of $f$, and $\Ht^k(F_{\!f,0},\Q)$ is identified with $\Hc^k\varphi_f\Q_{h,Y}$.
\sk
Set $Y':=Y\times_{\Delta}\Delta'$ with $\rho_m:(\Delta',0)\to(\Delta,0)$ a totally ramified $m$-fold covering. Let $\beta$ be the smallest positive rational number such that $\e(\beta)\,(:=e^{2\pi i\beta})$ is an eigenvalue of the Milnor monodromy of $f$. Assume
$$\h{$\frac{1}{m}$}\les\beta.
\leqno(4.8.2)$$
The following three conditions are then equivalent to each other:
\skn
(a) $(X,0)$ is a Du Bois singularity.
\skn
(b) $(Y',0)$ is a rational singularity.
\skn
(c) $f:Y\to\Delta$ is a cohomologically insignificant smoothing.
\skn
Condition~(c) means that $\Gr_F^0\Ht^k(F_{\!f,0},\C)=0$ ($\forall\,k$), see \cite{St3}. (This condition is invariant by the base change of $\Delta$.)
\sk
Set $h=f-z^m$ on $Y\times\C$ with $z$ the coordinate of $\C$ so that $Y'=h^{-1}(0)$. Then the above three conditions are respectively equivalent to
\skn\q\q
(a)$'\,\,\,\,\al_{f,1}\ges 1$.\q\q(b)$'\,\,\,\,\al_{h,1}>1$.\q\q(c)$'\,\,\,\,\al_{f,\mu_f}\les d_Y-1$.
\skn
Indeed, the first two equivalences follow from the arguments related to conditions (4.2.4), (4.3.8), and the last one from the above definition of spectrum, see (4.8.1). We have moreover the symmetry (see \cite{St2}):
$$\al_{f,i}+\al_{f,j}=d_Y\q\h{if}\q i+j=\mu_f+1,
\leqno(4.8.3)$$
together with the Thom-Sebastiani theorem as in \cite{ScSt}, \cite{Va}:
$${\rm Sp}(h)={\rm Sp}(f)\,{\rm Sp}(g),
\leqno(4.8.4)$$
where $g:=z^m$. Since ${\rm Sp}(g)=\msum_{k=1}^{m-1}\,t^{\,k/m}$ (see Remark~(3.3)(ii)), we then get
$$\al_{h,1}=\al_{f,1}+\h{$\frac{1}{m}$}.
\leqno(4.8.5)$$
So the equivalences between (a), (b), (c) follow.
\sk
In the case $\rho_m$ is associated with a semi-stable reduction, the above equivalences are a special case of \cite[Theorem 3.)]{St3} combined with \cite[Theorem 6.3]{Is} where an arbitrary smoothing of a normal (or Cohen-Macaulay) isolated singularity is treated. We take a projective compactification of $f$ as in \cite{Br} to apply \cite{Is}.


\begin{thebibliography}{CaMaScSh}
\bibitem[Al]{Al} Aluffi, P., Grothendieck classes and Chern classes of hyperplane arrangements, Int.\ Math.\ Res.\ Notices (2013), 1873--1900.
\bibitem[BaFuMa]{BFM} Baum, P., Fulton, W.\ and MacPherson, R., Riemann-Roch for singular varieties, Inst.\ Hautes Etudes Sci.\ Publ.\ Math.\ 45 (1975), 101--145.
\bibitem[BrScYo]{BSY} Brasselet, J.-P., Sch\"urmann, J. and Yokura, S., Hirzebruch classes and motivic Chern classes of singular spaces, Journal of Topology and Analysis 2, (2010), 1--55.
\bibitem[Br]{Br} Brieskorn, E., Die Monodromie der isolierten Singularit\"aten von Hyperfl\"achen, Manuscripta Math., 2 (1970), 103--161.
\bibitem[Bu]{Bu} Budur, N., On Hodge spectrum and multiplier ideals, Math.\ Ann.\ 327 (2003), 257--270.
\bibitem[BuSa1]{BS1} Budur, N.\ and Saito, M., Multiplier ideals, $V$-filtration, and spectrum, J.\ Alg.\ Geom.\ 14 (2005), 269--282.
\bibitem[BuSa2]{BS2} Budur, N.\ and Saito, M., Jumping coefficients and spectrum of a hyperplane arrangement, Math.\ Ann.\ 347 (2010), 545--579.
\bibitem[CaMaScSh]{CMSS} Cappell, S.E., Maxim, L., Sch\"urmann, J.\ and Shaneson, J.L., Characteristic classes of complex hypersurfaces, Adv.\ Math.\ 225 (2010), 2616--2647.
\bibitem[De1]{De1} Deligne, P., Th\'eorie de Hodge II, Publ.\ Math.\ IHES, 40 (1971), 5--58.
\bibitem[De2]{De2} Deligne, P., Le formalisme des cycles \'evanescents, in SGA7 XIII and XIV, Lect. Notes in Math. 340, Springer, Berlin, 1973, 82--115 and 116--164.
\bibitem[DeLo]{DL} Denef, J.\ and Loeser, F.: Motivic exponential integrals and a motivic Thom-Sebastiani theorem, Duke Math.\ J.\ 99 (1999), 285--309.
\bibitem[DiMaSaTo]{DMST} Dimca, A., Maisonobe, Ph., Saito, M.\ and Torrelli, T., Multiplier ideals, $V$-filtrations and transversal sections, Math.\ Ann.\ 336 (2006), 901--924.
\bibitem[Fu]{Fu} Fulton, F., Intersection Theory, Springer, Berlin, 1984.
\bibitem[FJ]{FJ} Fulton, W.\ and Johnson, K., Canonical classes on singular varieties, Manuscripta Math.\ 32 (1980), 381--389.
\bibitem[GeLoMe]{GLM} Guibert, G., Loeser, F., Merle, M., Iterated vanishing cycles, convolution, and a motivic analogue of a conjecture of Steenbrink, Duke Math.\ J.\ 132 (2006), 409--457.
\bibitem[Ha]{Ha} Hartshorne, R., Algebraic Geometry, Springer, New York, 1977.
\bibitem[Hi]{Hi} Hirzebruch, F., Topological methods in algebraic geometry, Springer, Berlin, 1966.
\bibitem[Is]{Is} Ishii, S., On isolated Gorenstein singularities, Math.\ Ann.\ 270 (1985), 541--554.
\bibitem[Ka]{Ka} Kashiwara, M., Vanishing cycle sheaves and holonomic systems of differential equations, Lect.\ Notes in Math. 1016, Springer, Berlin, 1983, pp. 136--142.
\bibitem[KoSc]{KoSc} Kov\'acs, S.\ and Schwede, K., Hodge theory meets the minimal model program: a survey of log canonical and Du Bois singularities, MSRI Publ., 58, Cambridge Univ.\ Press, 2011, pp.~51--94.
\bibitem[La]{La} Lazarsfeld, R., Positivity in algebraic geometry II, Springer, Berlin, 2004.
\bibitem[Mac]{Mac} MacPherson, R.D., Chern classes for singular algebraic varieties, Ann.\ of Math.\ (2) 100 (1974), 423--432.
\bibitem[Mal]{Mal} Malgrange, B., Polyn\^ome de Bernstein-Sato et cohomologie \'evanescente, Ast\'erisque 101-102 (1983), 243--267.
\bibitem[MaSaSc1]{MSS1} Maxim, L., Saito, M.\ and Sch\"urmann, J., Hirzebruch-Milnor classes of complete intersections, Adv.\ in Math.\ 241 (2013) 220--245.
\bibitem[MaSaSc2]{MSS2} Maxim, L., Saito, M.\ and Sch\"urmann, J., Hirzebruch-Milnor classes and Steenbrink spectra of certain projective hypersurfaces (arXiv:1312.0392), to appear in the Proceedings of the Arbeits\-tagung 2013.
\bibitem[MaSaSc3]{MSS3} Maxim, L., Saito, M.\ and Sch\"urmann, J., Thom-Sebastiani theorems for filtered $\D$-modules and for multiplier ideals (arXiv:1610.07295).
\bibitem[Nav]{Nav} Navarro Aznar, V., On the Chern classes and the Euler characteristic for nonsingular complete intersections, Proc.\ Amer.\ Math.\ Soc.\ 78 (1980), 143--148.
\bibitem[OhYo]{OY} Ohmoto, T.\ and Yokura, S., Product formulas for the Milnor class, Bull.\ Polish Acad.\ Sci.\ Math.\ 48 (2000), 387--401.
\bibitem[OrSo]{OS} Orlik, P.\ and Solomon, L., Combinatorics and topology of
complements of hyperplanes, Inv.\ Math.\ 56 (1980), 167--189.
\bibitem[PaPr]{PP} Parusi\'nski A.\ and Pragacz, P., Characteristic classes of hypersurfaces and characteristic cycles, J.\ Alg.\ Geom.\ 10 (2001), 63--79.
\bibitem[Sa1]{exg} Saito, M., On the exponents and the geometric genus of an isolated hypersurface singularity, in Singularities, Part 2, Proc.\ Sympos.\ Pure Math.\ 40, AMS Providence RI, 1983, pp.~465--472.
\bibitem[Sa2]{mhp} Saito, M., Modules de Hodge polarisables, Publ.\ RIMS, Kyoto Univ.\ 24 (1988), 849--995.
\bibitem[Sa3]{dual} Saito, M., Duality for vanishing cycle functors. Publ.\ RIMS, Kyoto Univ.\ 25 (1989), 889--921.
\bibitem[Sa4]{mhm} Saito, M., Mixed Hodge modules, Publ.\ RIMS, Kyoto Univ.\ 26 (1990), 221--333.
\bibitem[Sa5]{rat} Saito, M., On $b$-function, spectrum and rational singularity, Math.\ Ann.\ 295 (1993), 51--74.
\bibitem[Sa6]{mic} Saito, M., On microlocal $b$-function, Bull.\ Soc.\ Math.\ France 122 (1994), 163--184.
\bibitem[Sa7]{ste} Saito, M., On Steenbrink's conjecture, Math.\ Ann.\ 289 (1991), 703--716.
\bibitem[Sa8]{fil} Saito, M., On the Hodge filtration of Hodge modules, Moscow Math.\ J.\ 9 (2009), 161--191.
\bibitem[Sa9]{ts} Saito, M., Thom-Sebastiani theorem for Hodge modules, preprint.
\bibitem[ScSt]{ScSt} Scherk, J.\ and Steenbrink, J.H.M., On the mixed Hodge structure on the cohomology of the Milnor fibre, Math.\ Ann.\ 271 (1985), 641--665.
\bibitem[Sch1]{Sch1} Sch\"urmann, J., A generalized Verdier-type Riemann-Roch theorem for Chern-Schwartz-MacPherson classes (arXiv:math/0202175).
\bibitem[Sch2]{Sch2} Sch\"urmann, J., Topology of singular spaces and constructible sheaves, Birkh\"auser Verlag, Basel, 2003.
\bibitem[Sch3]{Sch3} Sch\"urmann, J., Characteristic classes of mixed Hodge modules, in Topology of Stratified Spaces, MSRI Publications Vol.\ 58, Cambridge University Press (2011), 419--471 (arXiv:0907.0584).
\bibitem[Sch4]{Sch4} Sch\"urmann, J., Specialization of motivic Hodge-Chern classes (arXiv:0909.3478).
\bibitem[Sch5]{Sch5} Sch\"urmann, J., Nearby cycles and characteristic classes of singular spaces, in IRMA Lectures in Mathematics and Theoretical Physics, European Math.\ Soc., Vol. 20 (2012), pp.~181--205 (arXiv:1003.2343).
\bibitem[Sch6]{Sch6} Sch\"urmann, J., Chern classes and transversality for singular spaces (arXiv:1510.01986), to appear in Singularities in Geometry, Topology, Foliations and Dynamics, Trends in Mathematics 2017.
\bibitem[SGA6]{SGA6} Berthelot, P., Grothendieck A.\ and Illusie, L., Th\'eorie des intersections et th\'eor\`eme de Riemann-Roch, SGA 6, Lect.\ Notes in Math.\ 225, Springer, Berlin, 1971.
\bibitem[St1]{St1} Steenbrink, J.H.M., Limits of Hodge structures, Inv.\ Math.\ 31 (1975/76), 229--257.
\bibitem[St2]{St2} Steenbrink, J.H.M., Mixed Hodge structure on the vanishing cohomology, in Real and complex singularities (Proc.\ Ninth Nordic Summer School,Oslo, 1976), Sijthoff and Noordhoff, Alphen aan den Rijn, 1977, 525--563.
\bibitem[St3]{St3} Steenbrink, J.H.M., Mixed Hodge structures associated with isolated singularities, in Singularities, Part 2, Proc.\ Sympos.\ Pure Math.\ 40, AMS Providence RI, 1983, pp.~513--536.
\bibitem[Va]{Va} Varchenko, A. N., Asymptotic mixed Hodge structure in vanishing cohomologies, Math.\ USSR Izv.\ 18 (1982), 469--512.
\bibitem[Ve]{Ve} Verdier, J.-L., Le th\'eor\`eme de Riemann-Roch pour les intersections compl\`etes, Ast\'erisque 36--37 (1976), 189--228.
\bibitem[Yo1]{Yo1} Yokura, S., A generalized Grothendieck-Riemann-Roch theorem for Hirzebruch's $\chi_y$-charac\-teristic and $T_y$-characteristic, Publ.\ RIMS, Kyoto Univ.\ 30 (1994), 603--610.
\bibitem[Yo2]{Yo2} Yokura, S., On characteristic classes of complete intersections, in Algebraic geometry: Hirzebruch 70 (Warsaw, 1998), Contemp. Math., 241, Amer. Math. Soc., Providence, RI, 1999, 349--369.
\end{thebibliography}
\end{document}